\documentclass[a4paper,twosided]{amsart}
\usepackage{amsmath,amssymb,relsize,color,amsthm,graphicx,float,rotfloat,afterpage}

\definecolor{refkey}{gray}{.40} 
\definecolor{labelkey}{gray}{.20} 

\newcommand{\ForJournal}[1]{{}}
\newcommand{\ForArxiv}[1]{{#1}}

\newcommand{\arXiv}{\textsf{arXiv}}

\usepackage[hyperindex, pdfmark, pdftex,colorlinks=true,
                      pdfstartview=FitV,
                      linkcolor=blue,
                      citecolor=blue,
                      urlcolor=blue,
          ]{hyperref}
\DeclareMathOperator{\hgt}{ht}
\DeclareMathOperator{\UnitOp}{U}

\DeclareMathOperator{\sgn}{sgn}
\DeclareMathOperator{\sortop}{sort}
\DeclareMathOperator{\id}{Id}

\newcommand{\fracG}{\frac{G(k+1)^2}{G(2k+1)}}

\newtheorem{thm}{Theorem}
\newtheorem{lemma}[thm]{Lemma}

\newtheorem{prop}[thm]{Proposition}
\newtheorem*{fact}{Fact}

\newcommand{\Ones}[1]{{\left(\left[1^{#1}\right]\right)}}

\newcommand{\boldmu}{\boldsymbol{\mu}}
\newcommand{\boldrho}{\boldsymbol{\rho}}
\newcommand{\boldalpha}{\boldsymbol{\alpha}}
\newcommand{\boldbeta}{\boldsymbol{\beta}}

\newcommand{\total}{\substack{\sum}}
\newcommand{\vect}[1]{{\overline{\mathbf{#1}}}}
\newcommand{\sort}[1]{{\overrightarrow{\mathbf{#1}}}}
\newcommand{\Unit}[1]{{\UnitOp(#1)}}

\newcommand{\ud}{\text{d}}
\newcommand{\averageLeft}{\left<}
\newcommand{\averageRight}{\right>_\Unit{N}}
\newcommand{\average}[1]{\averageLeft #1 \averageRight}

\newcommand{\Nplus}{\mathbb{N}_+}
\newcommand{\Frob}[2]{\left\{
  \substack{#1\\{#2}}\right\}}
\title[Joint moments]{Joint moments of derivatives\\ of characteristic
  polynomials}
\dedicatory{Pour Annie \& Jean-Paul}
\email{paul-olivier.dehaye@merton.ox.ac.uk}
\author{Paul-Olivier Dehaye}
\date{\today}
\address{University of Oxford, Merton College, Merton Street, OX1 4JD, Oxford, United Kingdom}

\begin{document}
\begin{abstract}
We investigate the joint moments of the $2k$-th power of the
characteristic polynomial of random unitary matrices with the $2h$-th
power of the derivative of this same polynomial. We prove that for a
fixed $h$, the moments are given by rational functions of $k$, up to a
well-known factor that already arises when $h=0$. 

We fully describe the denominator in those rational functions (this
had already been done by Hughes experimentally), and define the
numerators through various formulas, mostly sums over partitions. 

We also use this to formulate conjectures on joint moments of the zeta
function and its derivatives, or even the same questions for the Hardy
function, if we use a ``real'' version of characteristic polynomials.

Our methods should easily be applicable to other similar problems, for
instance with higher derivatives of characteristic polynomials.

\ForJournal{More data and computer programs are available as expanded
  content.}
\ForArxiv{More data is available online, either on the author's
  web site or attached to the \LaTeX\, source of this \arXiv\, submission.}
\end{abstract}

\maketitle

\tableofcontents
\section{Introduction}
In Section~\ref{results}, we merely define what is meant by \emph{joint moments of characteristic polynomials} and state the results obtained in this paper. In Section~\ref{motivation}, we motivate these Random Matrix Theory results by Number Theory questions and explain the interest of joint moments in the context of the Riemann $\zeta$-function. 
In Section~\ref{techniques}, we discuss our techniques,
which are essentially in Representation Theory and Algebraic
Combinatorics. The organization of this paper is summarized in Section~\ref{organization}.
\subsection{Presentation of results}
\label{results}
We take for the characteristic polynomial of a
$N\times N$ unitary matrix $U$
\begin{eqnarray*}
Z_U(\theta) & := & \prod_{j=1}^N \left(1-e^{\mathfrak{i}(\theta_j-\theta)}\right),
\end{eqnarray*}
where the $\theta_j$s are the eigenangles of $U$.

We define 
\begin{eqnarray}
V_U(\theta) & := & e^{\mathfrak{i} N (\theta + \pi)/2} e^{-\mathfrak{i}
  \sum_{j=1}^N \theta_j/2} Z_U(\theta).
\label{VAppear}
\end{eqnarray}
It is easily checked that for real $\theta$, $V_U(\theta)$ is real and
$|V_U(\theta)|=|Z_U(\theta)|$. 

In this paper, we will investigate the averages (with respect to Haar measure)
\begin{eqnarray*}
\left|\mathcal{M}\right|_N(2k,r)&:=&\average{\left|Z_U(0)\right|^{2k} \left| \frac{Z_U'(0)}{Z_U(0)} \right|^{r}},\\
\left(\mathcal{M}\right)_N(2k,r)&:=&\average{\left|Z_U(0)\right|^{2k} \left( \frac{Z_U'(0)}{Z_U(0)} \right)^{r}},\\
\left|\mathcal{V}\right|_N(2k,r)&:=&\average{\left|V_U(0)\right|^{2k} \left|
    \frac{V_U'(0)}{V_U(0)} \right|^{r}} 
\end{eqnarray*}
and their asymptotics
\begin{eqnarray*}
\left|\mathcal{M}\right|(2k,r)&:=&\lim_{N\rightarrow \infty} \left|\mathcal{M}\right|_N(2k,r)/N^{k^2+r},\\
\left(\mathcal{M}\right)(2k,r)&:=&\lim_{N\rightarrow \infty} \left(\mathcal{M}\right)_N(2k,r)/N^{k^2+r},\\
\left|\mathcal{V}\right|(2k,r)&:=&\lim_{N\rightarrow \infty} \left|\mathcal{V}\right|_N(2k,r)/N^{k^2+r}.
\end{eqnarray*}
It is easy to show (by expanding the Haar measure explicitly) that the
averages at finite $N$ only make sense when $2k-r>-1$. For the
asymptotics, the normalization by $N^{k^2+r}$ is due to
Hughes~\cite{HughesJoint} (and proved in this paper anyways). 

This and related problems have been looked at by Conrey, Rubinstein and Snaith~\cite{CRS},
Hughes~\cite{HughesThesis,HughesJoint}, Hughes, Keating and O'Connell \cite{HughesFirst}, 
Forrester and Witte~\cite{ForresterWitte},
Mezzadri~\cite{Mezzadri2003}. However, much mystery remains, in
particular for the dependency in $r$ (when $r {\in}
\mathbb{R}\setminus \mathbb{N}$).  

While $r \in \mathbb{R}\setminus \mathbb{N}$ remains out of reach, we
offer in this paper an alternative approach that uncovers some of the structure in those averages:
\begin{thm}
\label{thmComplexM}
For $r \in \mathbb{N}$ and $k \in \mathbb{C}$, the moments
$(\mathcal{M})(2k,r)$ are essentially given by rational functions,
i.e. as meromorphic functions of $k$ we have
\begin{eqnarray}
(\mathcal{M})(2k,r) &=& \left(-\frac{\mathfrak{i}}{2}\right)^r \fracG
\frac{X_r(2k)}{Y_r(2k)}, \label{polComplex}
\end{eqnarray}
where $X_r$ and $ Y_r$ are even monic polynomials with integer
coefficients and with $\deg X_r = \deg Y_r$ and $G$ is the Barnes $G$-function~\cite[Appendix]{HughesFirst}.

Moreover
$$
Y_r(u) = \prod_{\substack{1 \le a \le r-1\\a \text{ odd}}} (u^2-a^2)^{\alpha_{a}(r)},
$$
with the $\alpha_{a}(\cdot)$ given by
$$
\alpha_a(r) =
\left\lfloor\frac{-a+\sqrt{a^2+4r}}{2}\right\rfloor .
$$
We derive from this a similar result (Theorem~\ref{otherMoments},
page~\pageref{otherMoments}) for  $|\mathcal{M}|(2k,2h)$ and
$|\mathcal{V}|(2k,2h)$ (for $h$ an integer). Finally, we have explicit
expressions for $(\mathcal{M})(2k,r)$ given in
Theorem~\ref{concise}, page~\pageref{concise} and
Theorem~\ref{ninthTheorem}, page~\pageref{ninthTheorem}  which allow
us to
compute the $X_r(u)$s, as given in Table~\ref{dataComplexM}, page
\pageref{dataComplexM}, and additional data (available in Section~\ref{data}).
\end{thm}
This structure had been guessed a few years ago by Chris Hughes based
on computational evidence \cite{HughesJoint}. The author is deeply
thankful to him for freely sharing and explaining all of his previous unpublished work.

This paper was initiated at Stanford University\footnote{with support
from FRG DMS-0354662} while the
author was finishing his Ph.D., mostly worked on at Merton College,
University of Oxford  and finalized at Institut des Hautes Etudes Scientifiques, Bures-sur-Yvette. The author wishes to thank his hosting institutions for their support as well as his Ph.D. adviser, Dan Bump, and Persi Diaconis, Masatoshi Noumi and Peter Neumann for helpful discussions.
\subsection{Motivation}
\label{motivation}
Ever since the works by Keating and Snaith~\cite{KS1,KS2},
the Riemann $\zeta$-function can be (conjecturally but quantitatively)
better understood through modelling by characteristic polynomials of
unitary matrices. The classical example concerns moments. Let 
\begin{eqnarray*}
g(k) & := & \fracG,\\
a(k) & := & \prod_{p \text{ prime}} \left(1 -\frac{1}{p}\right)^{k^2}
\sum_{m=0}^\infty \left(\frac{\Gamma(m+k)}{m!\,\Gamma(k)}\right)^2 p^{-m}.
\end{eqnarray*}
Then one can prove (fairly immediately, using the Selberg integral) that 
\begin{eqnarray}
|M|(2k,0) &=& g(k) ,\label{defg}
\end{eqnarray}
which according to the Keating-Snaith philosophy leads to the
following conjecture (for $k>-1/2$):
\begin{eqnarray}
\frac{1}{T}\int_0^T \left| \zeta(\frac{1}{2} + \mathfrak{i}
  t)\right|^{2k} \ud t \sim_T g(k) a(k) \left(\log \frac{T}{2 \pi}
\right)^{k^2}. \label{KS}
\end{eqnarray}
The main point is thus that $a(k)$ is obtained by looking at primes,
while $g(k)$ is guessed at from the Random Matrix side. 

Observe also that Equations~(\ref{defg}) and then (\ref{KS}) can be
analytically continued in $k$. 

Many of the authors cited above have now shown that this philosophy
should be extended to derivatives of characteristic polynomials.  

In particular, $|\mathcal{M}|(2k,r)$ should show up as the RMT factor of\footnote{It
  is a conjecture of Hall \cite{Hall2004} and Hughes \cite{HughesThesis}
  that this is the appropriate normalization with respect to $T$.}
$$
\mathcal{I}(2k,r) := \lim_{T \rightarrow \infty} \frac{1}{T \left(\log\frac{T}{2\pi}\right)^{k^2+r}} \int_0^T \left|
  \zeta(\frac{1}{2} + \mathfrak{i}t)\right|^{2k-r}
\left|\zeta'(\frac{1}{2} + \mathfrak{i}t)\right|^r \ud t.
$$
and similarly $|\mathcal{V}|(2k,r)$ is needed for 
$$
\mathcal{J}(2k,r) := \lim_{T \rightarrow \infty} \frac{1}{T \left(\log\frac{T}{2\pi}\right)^{k^2+r}} \int_0^T \left|
  \mathfrak{Z}(\frac{1}{2} + \mathfrak{i}t)\right|^{2k-r}
\left|\mathfrak{Z}'(\frac{1}{2} + \mathfrak{i}t)\right|^r \ud t,
$$
where $\mathfrak{Z}$ is Hardy's function (the relationship of
$\mathfrak{Z}$ to $\zeta$ is analogous to the relationship of $V_U$ to
$Z_U$, i.e.  when $t \in \mathbb{R}$, $\mathfrak{Z}(\frac12 + \mathfrak{i}t) \in \mathbb{R}$ and $\pm \mathfrak{Z}(\frac12 + \mathfrak{i}t)=\left|\zeta(\frac12 + \mathfrak{i}t)\right|$ ). More precisely, it is expected that 
$$
\mathcal{I}(2k,r) = a(k)  \,\, |\mathcal{M}|(2k,r) \text{ and } \mathcal{J}(2k,r) = a(k) \,\, |\mathcal{V}|(2k,r). 
$$
Thus, Theorems~\ref{thmComplexM} and~\ref{otherMoments}
give us a conjectural handle on the moments of $\zeta$ and $\mathfrak{Z}$. 

One can compute some small cases (for integer $k$ and $r$) and show
that they agree with previous Number Theory (proved) results. This had already
been done before and is repeated in Table~\ref{summary}.

\begin{table}
\caption{Summary of results on $\mathcal{I}(2k,2h)$ and
  $\mathcal{J}(2k,2h)$, when $h \ne 0$. The values for
  $|\mathcal{M}|(2k,2h)$ and $|\mathcal{V}|(2k,2h)$ are given as
  obtained from Theorem~\ref{otherMoments}. Observe that the fifth
  column equals the product of the  third and fourth. The ``source''
  column refers to the paper where the result in the fifth column was
  published. }
\label{summary}
\begin{tabular}{c|c|c|c|c|c}
\hline
$k$ & $h$ & $a(k)$ &
$|\mathcal{M}|(2k,2h)$& $\mathcal{I}(2k,2h)$ & source  \\
\hline
&&&&&\\
1&1&1&$\left(\frac{1}{2^2}\right)
\frac{1^2}{1}\frac{(2 \times 1)^2}{(2^2-1^2)^1}$&$\frac{1}{3}$&\cite{Ingham}\\
&&&&&\\
2&1&$\frac{6}{\pi^2}$&$\left(\frac{1}{2^2}\right)\frac{1^2}{12}\frac{(2\times
  2)^2}{(4^2-1^2)^1}$&$\frac{2}{15 \pi^2}$&\cite{ConreyFourth}\\
&&&&&\\
2&2&$\frac{6}{\pi^2}$&$\left(\frac{1}{2^4}\right)
\frac{1^2}{12}\frac{(2\times 2)^4 - 8(2\times 2)^2 - 6}{(4^2-1^2)^1(4^2-3^2)^1}$&$\frac{61}{1680
  \pi^2}$&\cite{ConreyFourth}\\
&&&&&\\
\hline
\multicolumn{6}{c}{}\\
\hline
$k$ & $h$ & $a(k)$&$|\mathcal{V}|(2k,2h)$ & $\mathcal{J}(2k,2h)$ & source \\
\hline
&&&&&\\
1&1&1&$\left(\frac{2!}{1!2^3}\right)
\frac{1^2}{1}\frac{1}{(2^2-1^2)^1}$&$\frac{1}{12}$&\cite{Ingham}, see
also \cite{HughesJoint}\\
&&&&&\\
2&1&$\frac{6}{\pi^2}$&$\left(\frac{2!}{1!2^3}\right)\frac{1^2}{12}\frac{1}{(4^2-1^2)^1}$&$\frac{1}{120
  \pi^2}$&\cite{ConreyFourth}, see also \cite{HallWirtinger}\\
&&&&&\\
2&2&$\frac{6}{\pi^2}$&$\left(\frac{4!}{2!2^6}\right)
\frac{1^2}{12}\frac{1}{(4^2-1^2)^1(4^2-3^2)^1}$&$\frac{1}{1120
  \pi^2}$&\cite{ConreyFourth}, see also \cite{HallWirtinger}\\
&&&&&\\
\hline
\end{tabular}
\end{table}

However, while Keating and Snaith obtained a full conjecture for $\mathcal{I}(2k,0)$
and $\mathcal{J}(2k,0)$ by computing $|\mathcal{M}|(2k,0)$ and $|\mathcal{V}|(2k,0)$, for the case
of joint moments this goal remains  elusive. All the available formulas
for $|\mathcal{M}|(2k,r)$ or $|\mathcal{V}|(2k,r)$ are rather inadequate. In particular,
those formulas are limited to $r:=2h$ ($h$ an integer), they are hard
to compute for large values of $k$ and $h$, they obscure some of the
structure in the results, and finally they cannot be
analytically continued in~$h$. 

The analytic continuation would be important, because Conrey and Ghosh
have proved in \cite{ConreyGhosh} that (on RH) $$\mathcal{J}(2,1) =
\frac{e^2-5}{4 \pi}$$ and hence 
effectively conjectured\footnote{This is completely backwards from the
  usual flow  of conjectures \emph{from} Random Matrix Theory
  \emph{to} Number Theory, and possibly an unique instance of a
  reversal of this type.}
$$|\mathcal{V}|(2,1) = \frac{e^2-5}{4 \pi} $$ as well since $a(1)=1$. In order
to get this, we would need to have a sufficiently nice formula for
$|\mathcal{V}|(2k,2h)$ that would allow for analytic continuation in $h$. We
have simply been unable to do this but have no doubt that our results
should be helpful for that goal (see the connection with Noumi's work below).

On the other hand, the formulas obtained in Theorem~\ref{concise},
page~\pageref{concise} allow for much more effective computation than
possible before, and we can compute longer tables for the different moments
(see Section~\ref{data}). 

This numerical data is useful as well, as Hall has devised (around
2002) a method
that uses $\mathcal{J}(2k,2h)$ for all $0\le h\le k$ to produce a lower bound $\Lambda(k)$ on 
$$
\Lambda := \limsup_{n\rightarrow \infty} \frac{t_{n+1}-t_n}{\frac{2\pi}{\log t_n}},
$$
where the $t_n$ is the $n$-th positive real zero of
$\zeta(1/2+\mathfrak{i}t)$. It is probably good to insist that this
method does not depend on the Riemann Hypothesis, but only on
 values for moments! At the time of
writing \cite{Hall2004}, Hall only had the information he needed for
$k$ up to 2 (conjecturally, up to $6$). In Section~\ref{data}, we
present our conjectural data for $\mathcal{J}(2k,2h)$ as a direct
function of $k$ for $h$ up to 15 
\ForJournal{(also available online \cite{DehayeData} for $h$ up to 30).}
\ForArxiv{(see \cite{DehayeData} or the source of this \arXiv\,submission for data up to $h=30$).}
 For a fixed $h$, various conjectural formulas are also given in
this paper for $\mathcal{J}(2k,2h)$ as a function of $k$. 
This, combined with Hall's method, should lead to more (conjectural) lower bounds on
$\Lambda$. It is widely believed that $\Lambda = \infty$ so
potentially we could also see if Hall's method has any hope to reach
that, assuming only information on the $\mathcal{J}(2k,2h)$, but not
on the Riemann Hypothesis. In other words, it would also inform us on
the relationship between moment conjectures, the Riemann
Hypothesis and the conjecture $\Lambda=\infty$. We leave this to a further paper. 

Finally, Noumi in his book \cite{NoumiBook} investigates the
relationship between Painlev\'e equations and expressions similar to
one of the expressions we obtain for $(\mathcal{M})(2k,r)$, in
Theorem~\ref{ninthTheorem}. Connections of this sort have been
uncovered before (see \cite{ForresterWitte,ForresterWitte2} and works
of Borodin), but an approach through Noumi's ideas would be
original. One of our goals then would be to obtain analytic
continuation for  $(\mathcal{M})(2k,r)$ in $r$, which would again
allow to compute $|\mathcal{V}|(2,1)$. We also leave this for
further study.

Our techniques are quite disconnected from the original motivation, so
we discuss them separately.
\subsection{Techniques} 
\label{techniques}
As mentioned earlier, our techniques lie mostly in Representation
Theory and Algebraic Combinatorics. We look at the characteristic
polynomials or the derivatives as symmetric functions of the
eigenvalues of $U$, and express them in that way. We eventually
express those symmetric functions in the most natural basis to use,
the Schur functions. This basis is particularly
suitable since those functions are also (irreducible) characters of
unitary groups $\Unit{N}$. We find ourselves integrating irreducible
characters over their support (groups), which is very enviable!

In order to express all the different functions in this basis of Schur
functions, we use ideas present in a paper of Bump and
Gamburd~\cite{BumpGamburd} and the author's
thesis~\cite{DehayeThesis}. We will introduce those ideas as we need
them.

For a more thorough discussion of why a similar approach
should always be attempted and other examples of its applications,
please see the author's thesis and the results in \cite{DehayeWeyl}.

Once we have a concise expression for the various moments, we still
have to evaluate it. This will involve sums over partitions of values
of Schur functions. After reparametrizing those sums over the
Frobenius coordinates of the partitions, results of El-Samra and King
were immediately useful to obtain the Schur values, and results of
Borodin to handle the combinatorics of the sums. We then obtain a very
big sum for the moments (Theorem~\ref{thmDeterminantExact}), but that
can directly be evaluated on computer (and thus checked against small
$N$ results). After taking asymptotics, our
results start simplifying into Theorem~\ref{concise}, enough to prove
Theorem~\ref{thmComplexM} on the general shape of those moments. 
However, the best expression is probably obtained once we use Macdonald's ninth
variation of the Schur functions (Theorem~\ref{ninthTheorem}). 
\subsection{Organization of this paper}
\label{organization}
\begin{itemize}
 \item In Section~\ref{notation}, we introduce all the non-standard
   notation we will be using.
 \item In Section~\ref{relations}, we present the basic relations
   satisfied by the integrands 
$\left|Z_U(0)\right|^{2k} \left| \frac{Z_U'(0)}{Z_U(0)} \right|^{r}$,
$\left|Z_U(0)\right|^{2k} \left( \frac{Z_U'(0)}{Z_U(0)} \right)^{r}$
and 
$\left|V_U(0)\right|^{2k} \left|\frac{V_U'(0)}{V_U(0)} \right|^{r}$.
 \item In Section~\ref{symmetric}, we re-express the integrands as a
   sum in the Schur basis, in a
   way similar to Bump and Gamburd (via the Dual Cauchy Identity).
 \item In Section~\ref{computation}, we engage in a long computation to evaluate
   the result obtained in the previous section, mostly using results
   of El-Samra and King, and Borodin. 
 \item Section~\ref{together} merely serves to tie what has been done
   in Sections~\ref{symmetric} and~\ref{computation} into the proof of Theorem~\ref{thmComplexM}.
 \item In Section~\ref{data} we present the data we are now able to
   compute, and particularly discuss the position of the roots of
   $|\mathcal{V}|(2k,2h)$ in Section~\ref{position}.
 \item Section~\ref{alternative} describes two attempts to simplify
   our results further, one using Macdonald's ninth variation of the
   Schur functions, and the second imitating a proof of the Cauchy identity.
\end{itemize}
The bulk of this paper is contained in
Sections~\ref{symmetric} and \ref{computation}.
\section{Notation}
\label{notation}
We let $\Nplus$ be the set $\mathbb{N}\setminus{0}$. To avoid
confusion with the index $i$, we have $\mathfrak{i}^2=-1$.

We use
$\vect{v}$ for a generic \emph{vector} (of integers) $(v_1,\cdots, v_d)$, and
$\sort{v}$ for a sorted sequence of strictly decreasing integers
$v_1>v_2>\cdots>v_d$, which we call a \emph{Frobenius
  sequence}. Frobenius sequences are thus a special type of vectors.  

Sequences of weakly decreasing positive integers
amount to partitions, and we stick with classical notation for those,
i.e. $\lambda = (\lambda_1,\cdots,\lambda_{l(\lambda)})$, which
defines $l(\lambda)$. We also freely change our point of view to Young
tableaux when discussing partitions.  
Given a partition~$\lambda$ of $|\lambda|$, we denote its conjugate by
$\lambda^t$. Define two sequences $p_i := \lambda_i-i$, $q_i :=
\lambda^t_i -i$. They are 
strictly decreasing;  $\lambda_i$ and $\lambda^t_i$ are eventually 0, and hence $p_i=-i$ and $q_i = -i$ eventually. There exists $d$ such that $p_d
\ge 0 > p_{d+1}$ and $q_d \ge 0 > q_{d+1}$. We call $d$ the \emph{rank} of
$\lambda$.  We have that $\sort{p} =
(p_1,\cdots,p_d)$ and $\sort{q} = (q_1,\cdots,q_d)$ are Frobenius
sequences, and we call $\sort{p}$ and $\sort{q}$ the \emph{Frobenius coordinates} of the
partition $\lambda $. We write $\lambda =
\Frob{\sort{p}}{\sort{q}}$. 

Given $\vect{p}$, we define
$\sigma_{\vect{p}} \in \mathcal{S}_d$ such that
$\sortop(\vect{p}) := (p_{\sigma_{\vect{p}}(i)})$ is strictly
decreasing (and hence a Frobenius sequence). This is thus not defined
if $p_i=p_j$ while $i \ne j$. We set $\sgn(\vect{p}) := \sgn(\sigma_{\vect{p}})$, with the added convention
that $\sgn(\vect{p}) := 0$ if $\sigma_{\vect{p}}$ is not defined. 

If $\lambda$ and $\mu$ are partitions, $\lambda \cup \mu$ is the
partition obtained by taking the union of their parts. The partition
$\left<X^Y\right>$ has a $Y\times X$ rectangle for Young tableau.

We also use the notation $\left[1^R\right]$ for $R$ copies of 1, used
as argument to a (Schur) function.

\section{Basic relations among the integrands}
\label{relations}
We logarithmically differentiate Equation~(\ref{VAppear}) to obtain 
\begin{eqnarray}
\frac{V_U'(\theta)}{V_U(\theta)} &=& \frac{\mathfrak{i} N}{2} +
\frac{Z_U'(\theta)}{Z_U(\theta)} \label{relationSource}
\end{eqnarray}
and hence, when $\theta$ is real, 
\begin{eqnarray*}
\left|\frac{Z_U'(\theta)}{Z_U(\theta)}\right|^2 & =
&\left|\frac{V_U'(\theta)}{V_U(\theta)}\right|^2 + \frac{N^2}{4},\\
& =
&\left(\frac{V_U'(\theta)}{V_U(\theta)}\right)^2 + \frac{N^2}{4},\\
&=&\left(\frac{Z_U'(\theta)}{Z_U(\theta)}\right)^2 +\mathfrak{i}N\left(\frac{Z_U'(\theta)}{Z_U(\theta)}\right).
\end{eqnarray*}

These basic relations give
\begin{eqnarray}
\left|\mathcal{M}\right|_N(2k,2h) & = &\sum_{j=0}^h (\mathfrak{i}N)^{h-j}
\binom{h}{j} (\mathcal{M})_N(2k,h+j),\label{relationOne}\\
\left|\mathcal{M}\right|(2k,2h) & = &\sum_{j=0}^h \mathfrak{i}^{h-j}
\binom{h}{j} (\mathcal{M})(2k,h+j), \label{relationTwo}\\
\left|\mathcal{V}\right|_N(2k,2h) & = &\sum_{j=0}^{h} \binom{h}{j}
\left(\frac{-N^2}{4}\right)^{h-j} |\mathcal{M}|_N(2k,2j),\label{relationThree}\\
\left|\mathcal{V}\right|(2k,2h) & = &\sum_{j=0}^{h} \binom{h}{j}
\left(\frac{-1}{4}\right)^{h-j} |\mathcal{M}|(2k,2j)
\label{relationFour}.
\end{eqnarray}
These formulas are initially valid only when $h$ is a non-negative integer, but
the RHSs can be analytically continued by plugging in non-integer $h$
and extending the sum to infinity\footnote{Getting the correct
  analytic continuation can be tricky: The relation 
\begin{eqnarray}
\label{eqnNotContinuation}
\left|\mathcal{V}\right|(2k,2h) &=&
\sum_{j=0}^{2h}\binom{2h}{j} \left(\frac{\mathfrak{i}}{2}\right)^j (\mathcal{M})(2k,2h-j)
\end{eqnarray}
is also valid for integers $h$ but here the RHS does not analytically
continue in $h$ to the LHS, since we exploit
$\left|\frac{V_U'(\theta)}{V_U(\theta)}\right|^{2h}$ =
$\left(\frac{V_U'(\theta)}{V_U(\theta)}\right)^{2h}$ where it is
critical that $h$ be an integer.}. 
We see thus that computing $(\mathcal{M})_N(2k,r)$ would get us most
of the way to $|\mathcal{M}|_N(2k,2h)$ or $|\mathcal{V}|_N(2k,2h)$, and we now focus on the integrand  $\left|Z_U(0)\right|^{2k}\left(\frac{Z_U'(0)}{Z_U(0)}\right)^r$.
\section{Derivation into the Schur basis}
\label{symmetric}
The goal here is to follow ideas similar to Bump and
Gamburd's~\cite{BumpGamburd} in order to prove
Proposition~\ref{symexpres}, page~\pageref{symexpres}. One of their
main tools was the Dual Cauchy identity. We encourage the reader to
look at their first Proposition and Corollary for the unitary group,
since this is all we really exploit from that paper.

\begin{lemma}[Dual Cauchy identity]If $\{x_i\}$ and $\{y_j\}$ are
  finite sets of variables,
$$\prod_{i,j}{1+x_iy_j} = \sum_\lambda s_{\lambda^t}(x_i)
s_\lambda(y_j),$$
where the sum is over all partitions $\lambda$ and $s_\lambda$ is the
Schur polynomial.
\end{lemma}

We (they) apply this Lemma setting $\{x_j :=
e^{\mathfrak{i}\theta_j}:j \in [1,\cdots, N ]\}$ to be the set of
eigenvalues of $U$, and $\{y_j := 1 : j \in [1,\cdots, 2k]\}$. We chose
the notation $s_\lambda(U):=s_\lambda(e^{\mathfrak{i}\theta_1},\cdots,e^{\mathfrak{i}\theta_N})$. This
gives
\begin{eqnarray*}
\sum_\lambda \overline{s_{\lambda^t}(U)} s_\lambda\left([1^{2k}]\right) &=& \det(\id + \overline{U})^{2k} \\
&=&  \overline{\det(U)}^k \left|\det(\id + U)\right|^{2k} \\
&=& \overline{s_{\left<k^N\right>}(U)} \left|\det(\id + U)\right|^{2k}
\end{eqnarray*}
or (replacing $U$ by $-U$)
\begin{eqnarray*}
|Z_U(0)|^{2k}=\left|\det(\id - U)\right|^{2k} &=& (-1)^{kN}s_{\left<k^N\right>}(U) \sum_\lambda (-1)^{|\lambda|}\overline{s_{\lambda^t}(U)} s_\lambda\left([1^{2k}]\right).
\end{eqnarray*}

We can also re-express 
\begin{eqnarray}
\frac{Z_U'(0)}{Z_U(0)} &=& \sum_{j=1}^N
\frac{\mathfrak{i}e^{\mathfrak{i}\theta_j}}{1-e^{\mathfrak{i}
    \theta_j}} \notag\\
&=&\sum_{j=1}^N \mathfrak{i} \lim_{z {\rightarrow} 1^-}
\sum_{m=1}^\infty  z^m e^{\mathfrak{i} m \theta_j }\notag\\
&=&\mathfrak{i}\lim_{z {\rightarrow} 1^-}\sum_{m=1}^\infty  z^m p_m(U),\label{zLimit}
\end{eqnarray}
where $p_m(x_1,\cdots,x_N)$ is the $m$-th power sum
$x_1^m+\cdots+x_N^m$ and we have used the same convention as for
$s_\lambda(U)$ of inputting the eigenvalues. We will use the same
convention soon for the power sums $p_\lambda := \prod_i p_{\lambda_i}$.

In practice, we want the reader to just ignore the variable $z$ and
set it to 1. This will be justified a posteriori.

Putting everything together, we thus get for $\left|Z_U(0)\right|^{2k}\left(\frac{Z_U'(0)}{Z_U(0)}\right)^r $
\begin{eqnarray}
 (-1)^{kN}s_{\left<k^N\right>}(U) \left(\mathfrak{i}\sum_{m=1}^\infty
   p_m(U)\right)^r \sum_\lambda
 (-1)^{|\lambda|} s_\lambda\left([1^{2k}]\right)\overline{s_{\lambda^t}(U)}.
\label{hugeSym}
\end{eqnarray}

At this point, we will soon want to use the fact that the
$s_\lambda$ are characters of unitary groups.

Indeed, if $U \in \Unit{N}$ then when $l(\lambda)>N$, we have\footnote{This is a consequence of the fact that  $s_\lambda(x_1,\cdots,x_n)\equiv 0 $ if $l(\lambda)>n$.\label{length}} $s_\lambda(U)\equiv0$, but when $l(\lambda),l(\mu) \le N$, we have
$$
\average{s_\lambda(U) \overline{s_\mu(U)}}= \delta_{\lambda\mu},
$$
i.e. for large enough $N$, $s_\lambda$ is an irreducible character of $\Unit{N}$.
This orthogonality is obviously good for our purposes, but the only
obstacle is the need to express $s_{\left<k^N\right>}(U)
\left(\sum_{m=1}^\infty   p_m(U)\right)^r$ exclusively in terms of
Schur functions. This can be done and will require the
Murnaghan-Nakayama rule.

Let a ribbon be a connected Young skew-tableau not containing any $2\times
2$-block. If a ribbon contains $m$ blocks, it is called a $m$-ribbon. 
A first approximation to one version of the M-N rule says that $s_\lambda p_m$ is
given by a signed sum of $s_\mu$s, where $\mu$ runs through all
partitions obtained by adding a $m$-ribbon to $\lambda$.

If we average Expression~(\ref{hugeSym}) over $\Unit{N}$, we could
thus see $\lambda$ as running through all partitions obtained by
adding $r$ ribbons to the rectangle $\left<N^k\right>$ (this uses the
fact that this lax version of the M-N rule is invariant under
transpositions, since we have yet to discuss the signs).
There are more conditions however. We also need $l(\lambda^t)\le N$
(since otherwise $s_{\lambda^t}(U)\equiv 0$, as in footnote
\ref{length}), and we need  $l(\lambda)\le 2k$ (since otherwise
$s_\lambda\Ones{2k}=0$, again just as in footnote \ref{length}).   
In other words, $\lambda$ contains $\left<N^k\right>$ but is contained
in $\left<N^{2k}\right>$. There are only finitely many (ways to obtain) such
partitions, which will make the sum over $\lambda$s finite, and thus
only finitely many sets of lengths of the $r$ ribbons will
contribute. This justifies a posteriori setting $z$ to 1 in
Equation~(\ref{zLimit}), but \emph{only when we can apply the dominated
  convergence theorem}. This will only occur if we know of a bound on
the integrand independent of $z$ that is itself integrable. We can
pick $\left|Z_U(0)\right|^{2k} \left| \frac{Z_U'(0)}{Z_U(0)}
\right|^{r}$ whenever this is integrable, i.e. only when $2k-r>-1$. 

We now state a more precise version of the M-N rule.
\begin{thm}[Murnaghan-Nakayama]
\label{MNprecise}
Let $\lambda$ be a partition and $\vect{\boldrho}$ be a vector with
$|\lambda| = \sum_i \boldrho_i$. If
$\chi^\lambda_\vect{\boldrho}$ is the value of the irreducible character of
$\mathcal{S}_{|\lambda|}$ associated to $\lambda$ on the conjugacy
class of cycle-type $\sortop(\vect{\boldrho})$, then
\begin{eqnarray}
p_{\vect{\boldrho}} &=& \sum_\lambda \chi^\lambda_{\vect{\boldrho}} s_\lambda\label{chartable}
\end{eqnarray}
and (more importantly)
\begin{eqnarray}
\chi^\lambda_\vect{\boldrho} &=& \sum_S (-1)^{\hgt(S)}\label{fullMN}
\end{eqnarray}
summed over all sequences of partitions
$S=(\lambda^{(0)},\lambda^{(1)},\cdots,\lambda^{(r)})$ such that
$r:=l(\lambda)$,
$0=\lambda^{(0)}\subset\lambda^{(1)}\subset\cdots\subset\lambda^{(r)}
=\lambda$, and such that each $\lambda^{(i)}-\lambda^{(i-1)}$ is a
ribbon of length $\boldrho_i$, and $\hgt(S)=\sum_i\hgt(\lambda^{(i)}-\lambda^{(i-1)})$.
\end{thm}
We have not defined the height $\hgt$ of a ribbon, but rather than
doing so or  detailing the computation here, we only expose the idea: 
Equation~(\ref{chartable}) tells us that $(\sum_m p_m)^r$ can be
computed using the character values of symmetric groups, which can be
evaluated by summing over sequences of partitions
$(\lambda^{(0)},\cdots,\lambda^{(r)})$. For each such sequence, the
sequence $(\tilde{\lambda}^{(0)},\cdots,\tilde{\lambda}^{(r)})$, with
$\tilde{\lambda}^{(i)}:=\left<N^k\right>\cup\lambda^{(i)}$, would be
associated with the combinatorics of the expansion of the product 
in Expression~(\ref{hugeSym}). Indeed the combinatorics of ribbon is
unchanged under translations (down by $k$) as long as the partitions are
kept within a rectangle (actually, a horizontally bounded region). 

If the computation is explicitly carried out, we get the following
result:
\begin{prop}
\label{symexpres}
If $2k-r>-1$, we have
\begin{eqnarray}
(\mathcal{M})_N(2k,r) &=& (-\mathfrak{i})^r
\sum_{\vect{\boldmu} \in \Nplus^r} \, \,
\sum_{\substack{\lambda \text{ within} \\ k \times N}}
\chi^\lambda_\vect{\boldmu} \,\, s_{\left<N^k\right> \cup \lambda}\Ones{2k},\label{RealMomentSeries}
\end{eqnarray}
with the understanding that $\chi^\lambda_\vect{\boldmu} = 0$ if
$|\lambda| \ne \sum_i \boldmu_i$.
\end{prop}
For this result, we have preferred to index all the partitions
containing $\left<N^k\right>$ but contained in $\left<N^{2k}\right>$
as $\left<N^k\right>\cup\lambda$, for $\lambda \subset \left<N^k\right>$.

We are now left with the task of evaluating the RHS in
Equation~(\ref{RealMomentSeries}), which will turn out to be a
tedious process.
\section{Main computation}
\label{computation}
We are left with two problems. The first one is due to the
characters of the symmetric group. Those are of course desperately hard to
evaluate directly and individually. We are helped here because we will actually only evaluate
something close to 
$$
\sum_{\substack{\vect{\boldmu} \in \Nplus^l}} \, \,
\chi^\lambda_\vect{\boldmu} 
$$
for given $\lambda$. This amounts to computing the sum of values of
the character $\chi^\lambda$ over permutations with $l$ cycles. The
second issue is evaluating 
$
s_{\left<N^k\right> \cup \lambda}\Ones{2k}
$. The author had previously used the Weyl Dimension Formula to do
this (see \cite{DehayeThesis}). A formula giving that
dimension in terms of the Frobenius coordinates of $\lambda$ is
probably better adapted for our purposes. 

In addition, both ``problems'' combine extremely well, in that both
expressions should involve a sign, which turns out to be the
same. 

We will then sum our terms over all partitions, expressed in Frobenius
coordinates. This amounts to summing over possible ranks ($1\le d$) and
then pairs of Frobenius sequences of length $d$.

\subsection{The value of the Schur function in Frobenius coordinates}
\subsubsection{Dimension formula in Frobenius coordinates}
El-Samra and King~\cite{SamraKing} use the notation $D_R \Frob{p}{q}$ for $s_{\Frob{\sort{p}}{\sort{q}}}\Ones{R} $.

Assume $\Frob{\sort{p}}{\sort{q}}$ has $d$ Frobenius coordinates. They prove that 
\begin{multline}
s_{\Frob{\sort{p}}{\sort{q}}} \Ones{R} 
= 
\left|
  \frac{(R+p_i)!}{(R-q_j-1)!p_i!q_j!(p_i+q_j+1)}\right|_{d\times d}
 \label{SamraKingFormula}\\
 =  \prod_{i=1}^d\frac{(R+p_i)!}{(R-q_i-1)!p_i!q_i!} \prod_{1 \le i < j \le d}(p_i-p_j)(q_i-q_j) 
\prod_{i,j=1}^d \frac{1}{p_i+q_j+1}
\end{multline}
where the first expression is also known as the reduced determinantal
form (cf. Foulkes \cite{Foulkes}, as cited in \cite{SamraKing}).

It is a consequence of Cauchy's Lemma that the two expressions in Formula~(\ref{SamraKingFormula}) are
equivalent:
\begin{lemma}[Cauchy]
\label{CauchysLemma}
\begin{eqnarray}
\left|
  \frac{1}{p_i+q_j+1}\right|_{d\times d}
=\prod_{1 \le i < j \le d}(p_i-p_j)(q_i-q_j) 
\prod_{i,j=1}^d \frac{1}{p_i+q_j+1}.  \label{CauchysEqn}
\end{eqnarray}
\end{lemma}

Observe that Formula~(\ref{SamraKingFormula}) is positive (as it
should, given that it is also a dimension) because the $p_i$ and
$q_i$ are strictly decreasing. 

However, the RHS of Formula~(\ref{SamraKingFormula}) still makes sense if we
plug in unsorted vectors $\vect{p},\vect{q}$ (with even the possibility
of $i \ne j$ but $p_i=p_j$). Hence this can be used to define
$s_{\Frob{\vect{p}}{\vect{q}}} \Ones{R}$ as well, which is then
skew-symmetric in both the $p_i$s and the $q_i$s separately. 
This can be written \begin{eqnarray}
s_{\Frob{\vect{p}}{\vect{q}}}\Ones{R} & =& \sgn(\vect{p})
\sgn(\vect{q}) s_{\Frob{\sortop(\vect{p})}{\text{sort(}\vect{q})}}\Ones{R}. \label{setSkew1}
\end{eqnarray}
Observe that Formula~(\ref{setSkew1}) is still valid when $\sortop(\vect{p})$
or $\sortop(\vect{q})$ is not defined (this happens when two of the entries of
$\vect{p}$ or $\vect{q}$ are equal) thanks to $\sgn(\vect{p})
\sgn(\vect{q})=0$ (see conventions in Section~\ref{notation})!

Finally, it is helpful to remark that Formula~(\ref{SamraKingFormula}) for
$s_{\Frob{\vect{p}}{\vect{q}}} \Ones{R}$ can be seen as a product indexed by the sets $\vect{p} \cup \vect{q}$ and pairs in the set
$\left(\vect{p}\times \vect{p}\right) \cup \left(\vect{q} \times\vect{q}\right) \cup \left(\vect{p}\times\vect{q}\right)$. 

\subsubsection{Evaluation of  $s_{\left<N^k\right> \cup \lambda}\Ones{2k}
$}
We take $\lambda = \Frob{\sort{p}}{\sort{q}}$ to have $d$ Frobenius coordinates.

In total analogy with Equation~(\ref{setSkew1}), we first extend the definition of $s_{\left<N^k\right> \cup \lambda}$ and set
$$
s_{\left<N^k\right> \cup \Frob{\vect{p}}{\vect{q}}}:= \sgn(\vect{p})
\sgn(\vect{q}) s_{\left<N^k\right> \cup
  \Frob{\sortop(\vect{p})}{\sortop(\vect{q})}},
$$
with the understanding (as before) that the value of the RHS should be 0 if
$p_i=p_j$ (resp.~$q_i=q_j$) for  $i \ne j$. Again, this is
skew-symmetric in the $p_i$s and separately in the $q_i$s.

We have the following Lemma
\begin{lemma}
Let $\vect{p}$, $\vect{q}$ be vectors with $d$ coordinates. Then
\begin{multline}
s_{\left<N^k\right> \cup \Frob{\vect{p}}{\vect{q}}}\Ones{2k} =\\
s_{\left<N^k\right>}\Ones{2k} 
\left(\prod_{i=1}^d\frac{(N-p_i)^{(k)}(k-q_i)^{(k)}}{(p_i+k+1)^{(k)}(N+q_i+1)^{(k)}}
\right)  s_{\Frob{\vect{p}}{\vect{q}}}\Ones{2k}. \label{SamraKingWithRectangle}
\end{multline}
\end{lemma}
\begin{proof}
By skew-symmetry, we really only have to check this for $\Frob{\sort{p}}{\sort{q}}$. If we want to use Formula~(\ref{SamraKingFormula}), we should look
at the Frobenius coordinates of $\left<N^k\right> \cup \lambda$. This
would be rather unpleasant (particularly because the number of
Frobenius coordinates would change for fixed $N$ and $k$ according to the $\lambda$
considered). 

Let us look instead at:
\begin{eqnarray*}
\sort{x} &:=& \left( N+k-1,\cdots,N\right),\\
\sort{y} &:=& \left( 2k-1,\cdots,k\right),\\
\sort{\boldalpha} &:=& \sort{x} \cup \sort{p}\quad \text{(sorted)}\\
\text{and }\sort{\boldbeta} &:=& \sort{y} \cup \sort{q} \quad \text{(sorted)}.
\end{eqnarray*}
Then $\sort{\boldalpha}$ and $\sort{\boldbeta}$ are strictly decreasing, so
those are Frobenius coordinates. The partition corresponding to those
coordinates is  obtained geometrically by sticking a $\left<k^{2k}\right>$
block to the left of $\left<N^k\right> \cup \lambda$, or equivalently
to shifting $\left<N^k\right> \cup \lambda$ by $k$ spots to the right,
while considering $\lambda = (\lambda_1,\cdots,\lambda_k)$ to have
exactly $k$ parts (with some possibly empty). 

Because of this, we have (as in \cite[page 6]{BumpGamburd}): $$s_{\Frob{\sort{\boldalpha}}{\sort{\boldbeta}}}\Ones{2k} =
e_{2k}^k\Ones{2k} s_{\left<N^k\right> \cup \lambda}\Ones{2k}
 = s_{\left<N^k\right> \cup \lambda} \Ones{2k}
.$$

Additionally, $\Frob{\sort{x}}{\sort{y}}$ are the Frobenius
coordinates of 
$\left<(N+k)^k\right> \cup \left<k^k\right>$. Hence, for the same
reason as above, we have:
\begin{multline*}
s_{\Frob{\sort{x}}{\sort{y}}}\Ones{2k}
=  s_{\left<(N+k)^k\right> \cup \left<k^k\right>}\Ones{2k} \\ =
e_{2k}^k\Ones{2k} s_{\left<N^k\right>}\Ones{2k} = s_{\left<N^k\right>}\Ones{2k}.
\end{multline*}

When evaluating the product described in Equation~(\ref{SamraKingFormula}) using the
$\sort{\boldalpha}$ and $\sort{\boldbeta}$ coordinates, we have
a big product taken over the sets $\sort{\boldalpha},
\sort{\boldbeta}, \sort{\boldalpha}\times
\sort{\boldalpha},\sort{\boldbeta}\times\sort{\boldbeta}$ and $\sort{\boldalpha}\times
\sort{\boldbeta}$. We expand those index sets using $\sort{\boldalpha}
= \sort{x} \cup \sort{p}$ and $\sort{\boldbeta} = \sort{y} \cup
\sort{q}$. 

One can see that the products indexed by $\sort{p}$, $\sort{q}$,
$\sort{p}\times\sort{p}$, $\sort{q}\times\sort{q}$ and
$\sort{p}\times\sort{q}$ together give $s_{\Frob{\sort{p}}{\sort{q}}}
\Ones{2k}=s_\lambda \Ones{2k}$. 

Similarly, the products indexed by $\sort{x}$, $\sort{y}$,
$\sort{x}\times\sort{x}$, $\sort{y}\times\sort{y}$ and
$\sort{x}\times\sort{y}$ give $s_{\Frob{\sort{x}}{\sort{y}}}
\Ones{2k}=s_{\left<N^k\right>} \Ones{2k}$.

We are left with only ``cross-products'' to evaluate, for the index
sets $\sort{x}\times\sort{p}$,
$\sort{x}\times\sort{q}$,$\sort{y}\times\sort{p}$ and
$\sort{y}\times\sort{q}$. The definitions of $\sort{x}$ and $\sort{y}$
now give the result.
\end{proof}
\subsection{Sums of characters over conjugacy classes with same number  of cycles}
Assume $f\left(\Frob{\vect{p}}{\vect{q}}\right)$ is a function of pairs of
vectors of the same length (say $d$). One can set
$f(\lambda):=f\left(\Frob{\sort{p}}{\sort{q}}\right)$, where $\lambda=\Frob{\sort{p}}{\sort{q}}$.

The goal in this section is to evaluate sums of characters of the
general form
$$
\sum_{\vect{\boldmu} \in \Nplus^l} 
\chi^\lambda_{\vect{\boldmu}} f(\lambda).
$$
We will eventually take $f(\lambda) = s_{\left<N^k\right>\cup\lambda}\Ones{2k}$ but there is no reason to limit ourselves in that way for a while.

We rely on a few results of Borodin that give
a slightly different version of the Murnaghan-Nakayama rule.
\subsubsection{Definitions}
This is based on  \cite[around page~15]{Borodin2000} and \cite[around
page~6]{Borodin1998}. 
The relevant definitions (not included here) are \emph{fragment}, the different
\emph{block} types, the \emph{filling numbers}, \emph{filled}
\emph{structure}, \emph{sign of a structure}. 

Theorem~\ref{thmStructures} is almost in Borodin's work, and his
definitions are used in Proposition~\ref{propStructuresCount}. Both of
those results are used for Theorem~\ref{thmDeterminantExact}, which
can be read without looking at Borodin's papers. 

However, the first condition to have a fragment needs clarification in both
papers:
change 
\begin{quote}
{(1) there is exactly one hook block that precedes the
  others} 
\end{quote}
to 
\begin{quote}
{(1) there is exactly one hook block in each
  fragment. That hook block precedes any other block in the fragment}
\end{quote}
 
We also would like to correct a statement in \cite{Borodin2000}, in that linear
horizontal or vertical blocks are \emph{positive}, not just
\emph{non-negative} integers (in agreement with the other cited paper of
Borodin~\cite{Borodin1998}).

We can highlight one of the definitions: Any filled structure $T$ with
$d$ fragments produces a set of pairs 
$$
\left\{ (p_1,q_1), \cdots , (p_d,q_d)\right\}
$$
which consists of the filling $p$- and $q$- numbers of the fragments.

The sign of $T$ is defined as follows:
$$
\sgn(T) = \sgn(\vect{p}) \sgn(\vect{q}) (-1)^{\sum q_i + v(T)}, 
$$
where, as a reminder, the $\sgn$ inside the formula is 0 if $p_i=p_j$
(resp.~$q_i=q_j$) for
$i \ne j$.
\subsubsection{Simplified Murnaghan-Nakayama rule}
Although we haven't defined anything, we state Proposition 4.3, taken
from the first paper of Borodin:
\begin{prop}
For any two partitions $\lambda$ and $\rho$ with $|\lambda|=|\rho|$, we
have 
$$
\chi^\lambda_\rho = \sum_T \sgn T,
$$
where the sum is taken over all filled structures of cardinality $\rho
= (\rho_1,\cdots,\rho_l)$ such that the sequences $(p_1,\cdots,p_d)
$ and $(q_1,\cdots,q_d)$ of filling $p$-numbers and $q$-numbers of the structure $T$ coincide, up to a
permutation, with the Frobenius $p$-coordinates and
$q$-coordinates of the partition $\lambda$ (i.e. $\lambda = \Frob{\sortop(\vect{p})}{\sortop(\vect{q})}$).
\end{prop}
The proof of this Proposition is quite simple: going back to the
original presentation of the Murnaghan-Nakayama rule in terms of
hooks, Borodin analyzes what happens to Frobenius coordinates when
subtracting hooks/ribbons. Each such subtraction corresponds to a
block. There are three cases to distinguish: the hook/ribbon can be above or
below the ``Frobenius diagonal'' or even overlap it. Those cases
correspond respectively to linear horizontal blocks, linear vertical blocks, and
hook blocks. 

This Proposition, as stated in Borodin's work, is slightly restrictive: there is
no need for $\rho$ to be a partition. Let $\vect{\boldrho}=(\rho_1,
\cdots,\rho_l)$ to be a vector of positive integers and define (just
as in Theorem~\ref{MNprecise})
$\chi^\lambda_{\vect{\boldrho}}:=\chi^\lambda_{\sortop(\vect{\boldrho})}$. Then,
by summing over all vectors $\vect{\boldrho}$, we get
\begin{prop}
For any partition $\lambda$, we have 
$$
\sum_{ \vect{\boldrho} \in \Nplus^l} \chi^\lambda_\vect{\boldrho} = \sum_T \sgn T,
$$
where the sum is taken over all filled structures $T$ of $l$ blocks
and 
with filling $p$-numbers $(p_1,\cdots,p_d)$ and $q$-numbers $(q_1,\cdots,q_d)$ such that  $\lambda = \Frob{\sortop(\vect{p})}{\sortop(\vect{q})}$.
\label{Boring}
\end{prop}
Observe that $d$, the rank of $\lambda$, has to be less or equal to
$l$ in order to have a structure.

We now state the main Theorem we will use that is originated in
Borodin's work. 
\begin{thm}
Assume $f$ is skew-symmetric within its two vector entries (separately), i.e.
$f\left(\Frob{\sortop(\vect{p})}{\sortop(\vect{q})}\right)=
\sgn(\vect{p})\sgn(\vect{q})f\left(\Frob{\vect{p}}{\vect{q}}\right)
$. Then,
\begin{eqnarray*}
\sum_{\substack{\lambda \text{ within}\\ k\times N}}\,\sum_{{\vect{\boldrho}\in \Nplus^l}} \chi^\lambda_{\vect{\boldrho}} f(\lambda) &=&
\sum_{d=1}^l \,\, \sum_{\substack{\vect{p}  \in [0,N-1]^d\\\vect{q}
    \in [0,k-1]^d}}
f\left(\Frob{\vect{p}}{\vect{q}}\right)
\sum_{T(\vect{p},\vect{q})} (-1)^{\sum q_i +v(T)},
\end{eqnarray*}
\label{thmStructures}
where $T(\vect{p},\vect{q})$ goes trough all filled structures of $d$ fragments,
$l$ blocks, $v(T)$ vertical blocks with filling $p$-numbers
$(p_1,\cdots,p_d)$ and $q$-numbers $(q_1,\cdots,q_d)$.
\end{thm}
\begin{proof}
We start by summing Proposition~\ref{Boring} over $\lambda$s fitting
inside a $k \times N$ box:
\begin{eqnarray*}
\sum_{\substack{\lambda \text{ within}\\ k\times N}}\,\sum_{{\vect{\boldrho}\in \Nplus^l}} \chi^\lambda_{\vect{\boldrho}} f(\lambda) &=&
\sum_{\substack{\lambda \text{ within}\\ k\times N}}  \sum_{T(\vect{p},\vect{q})}
(-1)^{\sum q_i +v(T)} \sgn(\vect{p}) \sgn(\vect{q}) \,
f\left(\lambda\right)\\
&=&
\sum_{\substack{\lambda \text{ within}\\ k\times N}}  \sum_{T(\vect{p},\vect{q})}
(-1)^{\sum q_i +v(T)} \, f\left(\Frob{\vect{p}}{\vect{q}}\right),
\end{eqnarray*}
where the second sums in each RHS are taken over all filled
structures $T(\vect{p},\vect{q})$ of $l$ blocks and $d$ fragments
such that the sequences of filling $p$-numbers $(p_1,\cdots,p_d)$
and $q$-numbers $(q_1,\cdots,q_d)$ of the structure 
coincide, up to two permutations, with the sequences of
Frobenius $p$-coordinates and $q$-coordinates of the partition
$\lambda$ (i.e. $\lambda =
\Frob{\sortop{(\vect{p})}}{\sortop{(\vect{q})}}$). Note that $d$ changes with $\lambda$.

We then obtain the final result by seeing the double sum over
$\lambda$ then permuted Frobenius coordinates of $\lambda$ as a sum over \emph{all} vectors of appropriate lengths. 

We should not be concerned about vectors having two identical coordinates
(say $p_i=p_j$), since the corresponding term in the RHS vanishes by
skew-symmetry of~$f$. 
\end{proof}
\subsubsection{Counting structures}
We now need to compute the sum
$$\sum_{T(\vect{p},\vect{q})} (-1)^{\sum q_i +v(T)}, $$
which is taken over the structures described above, i.e. for given
$l$, $d$, $\vect{p}$, $\vect{q}$, $v$. It would help to know how many
structures there are for each choices of those parameters. We prove the following Proposition.
\begin{prop}
There are exactly \label{propStructuresCount}
\begin{multline}
\# T(l,d,\vect{p},\vect{q},v) = \\
\sum_{\substack{\vect{s},\vect{t} \in \mathbb{N}^d\\ \sum t_i =
    v\\d+\total s_i+t_i = l}}\quad
\left[\frac{(s_d+t_d+\cdots+s_1+t_1+d)!}{\prod s_i! \prod t_i!
  \prod_{i=1}^d (d+1-i+\sum_{j=i}^d s_j+t_j)}\right]\times\left[
 \prod_i^d \binom{p_i}{s_i} \binom{q_i}{t_i} \right]\label{StructuresCount}
\end{multline}
structures with $d$ fragments, $l$ blocks, filling numbers
$\vect{p} = (p_1,\cdots,p_d)$ and $\vect{q}=(q_1,\cdots,q_d)$ and  $v$ vertical blocks.
The indices in the sum $s_i$ (resp. $t_i$) count horizontal (resp. vertical)
blocks in the $i^\text{th}$ fragment. 
\end{prop} 
\begin{proof}
This is a purely combinatorial problem. Given the number of vertical
blocks on each fragment, we essentially have a partial order on blocks
that we want to extend to form a linear order (across fragments). Part
of the rules in the initial partial order say that the hook-block in
the $i\text{th}$ fragment precedes any other block in that fragment. We
then need to fill the structure (i.e. choose filling numbers for each
block).

We can reverse this process:
\begin{itemize}
\item We first choose the numbers of horizontal
and vertical blocks $s_i$ and $t_i$ on the $i^\text{th}$ fragment. We have the
conditions that $\sum t_i = v(T)$ and $d + \sum s_i + t_i=l$
(i.e. there are $l$ blocks in total, $d$ hook, $s_i$ horizontal in the
$i^\text{th}$ fragment and $t_i$ horizontal in the $i^\text{th}$ fragment).  
\item Starting from the $d^\text{th}$ fragment, we decide where to insert the horizontal
  and vertical blocks of the $i^\text{th}$ fragment in the partial order
  that is established so far on the set of fragments from the $i+1$st
  to the $d^\text{th}$ one.
\item We decide how to cut up the $i^\text{th}$ fragment into filled blocks,
  respecting the number of horizontal/vertical blocks decided upon earlier.
\end{itemize}

The equality in the statement is intended to reflect clearly the
layering described above: the sum corresponds to the first layer, while the
other two layers correspond to one square-bracketed factor each.

Observe that the relation $s_d+t_d+\cdots+s_1+t_1+d = l$ could be used
to simplify the numerator in this expression.

The only hard part is to derive for the second step
\begin{multline*}
\frac{(s_d+t_d+\cdots+s_1+t_1+d)!}{\prod s_i! \prod t_i!
  \prod_{i=1}^d (d+1-i+\sum_{j=i}^d s_j+t_j)} = \\
\frac{(s_d+t_d+\cdots+s_1+t_1+d)!}{\prod s_i! \prod t_i!
 \left(\substack{ (s_d+t_d+1)\times\\(s_d+t_d+s_{d-1}+t_{d-1}+2)\times\\
     \cdots \times\\(s_d+t_d+s_{d-1}+t_{d-1}+\cdots+s_1+t_1+d)}\right)}.
\end{multline*}
This is obtained by simplifying
$$
\prod_{i=0}^{d-1} \binom{i+\sum_{j=d-i}^d s_j+t_j}{s_{d-i}+t_{d-i}}\binom{s_{d-i}+t_{d-i}}{s_{d-i}},
$$
where the $i^\text{th}$ factor in the $\prod_{i=0}^{d-1}$-product counts the
number of ways of choosing the linear order on the blocks of the
$d-i^\text{th}$ fragment, knowing the linear order restricted on the blocks of
the fragments $d-i+1$ to $d$. 

The first binomial factors intersperses the set of blocks of the $d-i^\text{th}$
fragment among the blocks of fragments $d-i+1$ to $d$, while the
second factor decides which blocks are horizontal and which are
vertical.
\end{proof}
We wish to insist on the fact that the summand in
Equation~(\ref{StructuresCount}) is not symmetric in the $p_i$s or the $q_i$s, because the factor in the denominator $\prod_{i=1}^d
(d+1-i+\sum_{j=i}^d s_j+t_j)$ is not symmetric in the $s_j$s or the
$t_j$s. For instance, $s_d$ appears $d$ times while $s_1$ appears
only once.
\subsubsection{Sum of determinants}
We aim now to put together all the results obtained so far in this section, but we first need a quick lemma.
\begin{lemma}
\label{LemmaSimpler}
Let $\vect{s}$ and $\vect{t}$ be vectors of integers. Then,
\begin{multline}
\sum_{\sigma,\tau \in \mathcal{S}_d} 
\frac{(\sgn{\sigma}\sgn{\tau})}{\prod_{i=1}^d
 (d+1-i+\sum_{j=i}^d s_{\sigma(j)}+ t_{\tau(j)})} = \\
\prod_{1\le i < j \le d} (s_i-s_j)(t_i-t_j) \prod_{1 \le i,j \le d}
\frac{1}{1+s_i+t_j}.
\label{simpler}
\end{multline}

\end{lemma}
\begin{proof}
The proof proceeds as for the classical computation for the Vandermonde determinant:
  the LHS is skew-symmetric in $\vect{s}$ and $\vect{t}$ separately,
  has obvious poles as prescribed in the RHS (when
  $s_{i_0}+t_{j_0}=-1$), and the degrees in the RHS are
  appropriate. Up to a constant of proportionality, both sides are thus the same. This
  constant is shown to be 1 by looking at rates of decrease when $s_1$
  goes to infinity.
\end{proof}

\begin{prop}
\label{propPainful}
Assume $f$ is skew-symmetric within its two vector entries (separately), i.e.
$f\left(\Frob{\sortop(\vect{p})}{\sortop(\vect{q})}\right)=
\sgn(\vect{p})\sgn(\vect{q})f\left(\Frob{\vect{p}}{\vect{q}}\right)
$. Then,
\begin{multline*}
\sum_{\vect{\boldmu} \in \Nplus^l} \, \,
\sum_{\substack{\lambda \text{ within} \\ k \times N}}
\chi^\lambda_\vect{\boldmu} f(\lambda) 
= l!
\sum_{d=1}^l \,\,
 \sum_{\substack{\vect{p}  \in [0,N-1]^d\\\vect{q}
     \in [0,k-1]^d\\v}}f\left(\Frob{\vect{p}}{\vect{q}}\right)
  (-1)^{\sum q_i +v}  \\
\sum_{\substack{\sort{s},\sort{t} \in \mathbb{N}^d\\ \sum t_i =
    v\\d+\total s_i+t_i = l}}
\prod_i^d \binom{p_{i}}{s_{i}} \binom{q_{i}}{t_{i}}
\prod_i \frac{1}{s_i!t_i!}\prod_{1\le i < j \le d} (s_i-s_j)(t_i-t_j) \prod_{1 \le i,j \le d}
\frac{1}{1+s_i+t_j}.
\end{multline*}
\end{prop}
\begin{proof}
We go head first and combine Equation~(\ref{RealMomentSeries}) and
Theorem~\ref{thmStructures}:\begin{multline*}
\sum_{\vect{\boldmu} \in \Nplus^l} \, \,
\sum_{\substack{\lambda \text{ within} \\ k \times N}}
\chi^\lambda_\vect{\boldmu} f(\lambda) = l!
\sum_{d=1}^l \,\, \sum_{\substack{\vect{p}  \in [0,N-1]^d\\\vect{q}
    \in [0,k-1]^d\\v}}f\left(\Frob{\vect{p}}{\vect{q}}\right)
 (-1)^{\sum q_i +v}  \\
\sum_{\substack{\vect{s},\vect{t} \in \mathbb{N}^d\\ \sum t_i =
    v\\d+\total s_i+t_i = l}}\quad
\frac{\prod_i^d \binom{p_i}{s_i} \binom{q_i}{t_i} }{\prod s_i! \prod t_i!
  \prod_{i=1}^d (d+1-i+\sum_{j=i}^d s_j+t_j)}\\
= l!
\sum_{d=1}^l \,\, \sum_{\substack{\vect{p}  \in [0,N-1]^d\\\vect{q}
    \in [0,k-1]^d\\v}}f\left(\Frob{\sortop(\vect{p})}{\sortop(\vect{q})}\right)
 (-1)^{\sum q_i +v}  \\
\sum_{\substack{\vect{s},\vect{t} \in \mathbb{N}^d\\ \sum t_i =
    v\\d+\total s_i+t_i = l}}\quad
\frac{\sgn(\vect{p})\sgn(\vect{q})\prod_i^d \binom{p_i}{s_i} \binom{q_i}{t_i} }{\prod s_i! \prod t_i!
  \prod_{i=1}^d (d+1-i+\sum_{j=i}^d s_j+t_j)}\\
= l!
\sum_{d=1}^l \,\,
\sum_{\substack{\Frob{\sort{p}}{\sort{q}}\\\text{within }
k\times N \\\vspace{.02in}\\v}}f\left(\Frob{\sort{p}}{\sort{q}}\right)
 (-1)^{\sum q_i +v}  \\
\sum_{\substack{\sort{s},\sort{t} \in \mathbb{N}^d\\ \sum t_i =
    v\\d+\total s_i+t_i = l}} \sum_{\sigma, \tau \in \mathcal{S}_d}\quad
\frac{\left[\sum_{\pi,\theta \in \mathcal{S}_d}\sgn(\pi)\sgn(\theta)\prod_i^d \binom{p_{\pi(i)}}{s_{\sigma(i)}} \binom{q_{\theta(i)}}{t_{\tau(i)}} \right]}{\prod s_{\sigma(i)}! \prod t_{\tau(i)}!
  \prod_{i=1}^d (d+1-i+\sum_{j=i}^d
  s_{\sigma(j)}+t_{\tau(j)})}\end{multline*}
It is now crucial to observe that for fixed
$\sort{p},\sort{q},\sort{s},\sort{t}$, the sign of the numerator of the summands
(bracketed) will depend on the parity of $\sigma$ and $\tau$. Hence we
obtain
\begin{multline*}
= l!
\sum_{d=1}^l \,\,
\sum_{\substack{\Frob{\sort{p}}{\sort{q}}\\\text{within }
k\times N \\\vspace{.02in}\\v}}f\left(\Frob{\sort{p}}{\sort{q}}\right)
 (-1)^{\sum q_i +v}  \\
\sum_{\substack{\sort{s},\sort{t} \in \mathbb{N}^d\\ \sum t_i =
    v\\d+\total s_i+t_i = l}}\sum_{\pi,\theta \in \mathcal{S}_d}\sgn(\pi)\sgn(\theta)\prod_i^d \binom{p_{\pi(i)}}{s_{i}} \binom{q_{\theta(i)}}{t_{i}}\\ \sum_{\sigma, \tau \in \mathcal{S}_d}\quad
\frac{\sgn(\sigma)\sgn(\tau) }{\prod s_{\sigma(i)}! \prod t_{\tau(i)}!
  \prod_{i=1}^d (d+1-i+\sum_{j=i}^d s_{\sigma(j)}+t_{\tau(j)})}
\end{multline*}
\begin{multline*}
= l!
\sum_{d=1}^l \,\,
\sum_{\substack{\Frob{\sort{p}}{\sort{q}}\\\text{within }
k\times N \\\vspace{.02in}\\v}}f\left(\Frob{\sort{p}}{\sort{q}}\right)
 (-1)^{\sum q_i +v}  
\sum_{\substack{\sort{s},\sort{t} \in \mathbb{N}^d\\ \sum t_i =
    v\\d+\total s_i+t_i = l}}
\\ \sum_{\sigma, \tau \in \mathcal{S}_d}\quad
\frac{\sgn(\sigma)\sgn(\tau) \prod_i^d \binom{p_{i}}{s_{\sigma(i)}} \binom{q_{i}}{t_{\tau(i)}} }{\prod s_{\sigma(i)}! \prod t_{\tau(i)}!
  \prod_{i=1}^d (d+1-i+\sum_{j=i}^d s_{\sigma(j)}+t_{\tau(j)})}.
\end{multline*}
The last line is now perfectly cut for substitution using 
Lemma~\ref{LemmaSimpler}. After changing the range of summation on
$\Frob{\sort{p}}{\sort{q}}$ within $k\times N$ to $\vect{p}  \in
[0,N-1]^d$,$\vect{q} \in [0,k-1]^d$, we obtain the announced result.

Admittedly, this is not very enlightening. It is thus worth
highlighting what happens: the sums we are dealing with initially are sums over
partitions. By using Frobenius coordinates, and sorting the partitions
by their rank $d$, we are expressing the main sum into a sum over $d$
of multisums in $d$ variables. We thus now have sums over two sets of
$d$ strictly decreasing variables (the sets $\sort{p}$ and
$\sort{q}$) of different ways of
building up this partition (the data encoded in $\vect{s}$ and
$\vect{t}$). Using skew-symmetry, we can unsort the variables
$\sort{p}$ and $\sort{q}$ to $\vect{p}$ and $\vect{q}$ and decide
instead to sort the variables according to ``building blocks'',
i.e. switch from $\vect{s}$ and $\vect{t}$ to $\sort{s}$ and $\sort{t}$.
\end{proof}
\subsection{Putting everything together}
We combine all the information obtained so far, and simultaneously
clear the restriction $d+\sum s_i+t_i = l$ in
Formula~(\ref{StructuresCount}) by encoding all the moments at once
into an exponential generating function.
\begin{thm}
For a fixed $k \in \mathbb{N}$, the two series  
\begin{multline}
\sum_{r>0} (\mathcal{M})_N(2k,r) \frac{(\mathfrak{i}z)^r}{r!}
\text{\quad and \quad} 
s_{\left<N^k\right>}\Ones{2k} 
\sum_{d=1}^\infty 
\sum_{\substack{\sort{s},\sort{t} \in \mathbb{N}^d}}
\\ \left|\frac{ z^{1+s_i+t_j}
}{s_i!  t_j!(1+s_i+t_j)}\right|_{d\times d}
 \sum_{\substack{\vect{p}  \in [0,N-1]^d\\\vect{q}
    \in [0,k-1]^d}}
\left| k\frac{ \binom{p_i}{s_i} \binom{q_j}{t_j}\binom{k+p_i}{p_i} \binom{k-1}{q_j}
    (N-p_i)^{(k)}(-1)^{q_j}}{(N+q_j+1)^{(k)}(1+p_i+q_j)}
\right|_{d\times d}
\label{seriesDeterminant}
\end{multline}
\label{thmDeterminantExact}
have equal coefficients of $z^r$ for $r < 2k+1$. 
\end{thm}
\begin{proof}
A first necessary remark is that as a formal power series, the second series is
well-defined: the sum to obtain the $r^\text{th}$ coefficient in that series reduces to a finite sum (because $s_i \le p_i$ and $t_j \le q_j$).

We know from Equation~(\ref{setSkew1}) that $s_{\left<
    N^k\right> \cup \Frob{\vect{p}}{\vect{q}}}\Ones{2k}$ is
skew-symmetric in $\vect{p}$ and $\vect{q}$ (separately). Hence we can
combine Equations~(\ref{RealMomentSeries}), (\ref{SamraKingFormula})
and (\ref{SamraKingWithRectangle}) and
Proposition~\ref{propPainful} to obtain a huge sum. The main statement
then follows from recombinations of the main product into
determinants, using Cauchy's Lemma~\ref{CauchysLemma}.
\end{proof}
\textbf{Remarks on Theorem~\ref{thmDeterminantExact}}
\begin{itemize}
\item This is a hypergeometric multisum (at least for fixed $d$), when we expand
the determinants using Cauchy's Lemma. However, not even small $d$s
seem tractable on computer.
\item A definite advantage of this formula is that it can be tested at
finite $N$ (by expanding the integral defining $(M_N)(2k,r)$
symbolically using the Haar measure). This is helpful to confirm the
results obtained so far.
\item We wish to insist on the idea behind this Theorem: initially we
  had a combinatorial problem on structures (see Formula
  (\ref{StructuresCount}) that had no symmetry for
  its summands in the $s_i$s (resp. $t_i$s). We have exploited some
  skew-symmetry in the $a$s and $b$s in
  Formula~(\ref{SamraKingWithRectangle}) to change this. In
  particular, we have now switched from a sum over
  $\sort{p},\sort{q},\vect{s},\vect{t}$ to a sum over
  $\vect{p},\vect{q},\sort{s},\sort{t}$. We have also simplified the denominator in Formula~(\ref{StructuresCount}). 
\item As a consequence of the previous point, we can now assume that
  the $s_i$s are all different. The same is true for the $t_i$s.
\item This has useful consequences, especially for computational
  purposes.  It is interesting to compute a bound on $r$
such that partitions with $d$ fragments will have a non-zero
contribution to the final sum in $(\mathcal{M})_N(2k,r)$. We have $r\ge d + \sum s_i + t_i$, and the
$s_i$s (resp. $t_i$s) should be all different. We can take them to
be $0,1,\cdots,d-1$. We thus have $r\ge d +2 \frac{d(d-1)}{2} =
d^2$.
\end{itemize}

We now define 
\begin{multline}
\label{HNkst}
H^{N,k,s,t} := \\{s!t!}\sum_{\substack{p  \in [0,N-1]\\ q \in [0,k-1]}}\frac{k
  (N-p)^{(k)}(-1)^{q}}{(N+q+1)^{(k)}(1+p+q)}\binom{k+p}{p}
\binom{k-1}{q} \binom{p}{s} \binom{q}{t},
\end{multline}
where the RHS is taken to be similar to the entries in one of the
determinants in Equation~(\ref{seriesDeterminant}). 

I have not been able to obtain a much better expression for this with
\textsf{Mathematica}. Normally, the package \textsf{MultiSum} \cite{MultiSum} should be able to deal with multiple hypergeometric series, but this
particular one is too complicated. We will thus focus on an easier
problem from now on, the problem of asymptotics (i.e. we switch from
$(\mathcal{M})_N(2k,r)$ to $(\mathcal{M})(2k,r)$).


\subsection{Asymptotics}
We need to compute asymptotics for $H^{N,k,s,t}$ more precisely.
\begin{prop}
For a fixed integer $k \ge 1$, when $k>t$,
\label{hypergeometricH}
\begin{eqnarray}
H^{k,s,t} & := & \lim_{N\rightarrow \infty} \frac{H^{N,k,s,t}}{N^{1+s+t}}\label{identityOne}\\
 & =&  k \sum_{i=0}^{k-t-1} \frac{\Gamma(k+i) 
  \Gamma(s+i+t+1)}{\Gamma(i+1)\Gamma(k+s+t+i+2)}\label{identityTwo}\\
&=& \frac{1}{1+s+t}\frac{\prod_{i=k}^{2k-1}
  (i-t)}{\prod_{j=k+1}^{2k}(j+s)}\label{identityThree}\\
&=
\footnotemark
&
\frac{1}{1+s+t}\frac{\Gamma(2k-t)\Gamma(k+s+1)}{\Gamma(k-t)\Gamma(2k+s+1)}.
\label{identityFour}
\end{eqnarray}
\footnotetext{Observe that Expression~(\ref{identityFour}) is well
  defined, thanks to the bound $k>t$. }
\end{prop}
\begin{proof}
Define 
$$
\tilde{H}^{N,k,s,t}  := t! \sum_{\substack{p  \in [0,N-1]\\ q \in [0,k-1]}}\frac{k
  (N-p)^{k}(-1)^{q}}{(N+q+1)^{k}(1+p+q)}\frac{p^k}{k!}
\binom{k-1}{q} {p^s} \binom{q}{t},
$$
i.e. $H^{N,k,s,t}$ stripped of some of its terms of obviously lower
order in $p$, $N$ and $q$ combined. We do this because we want to
compute the leading order of $H^{N,k,s,t}$ and there will be lots of
cancellation due to the sum over $q$ (as showed by
Equation~(\ref{cancellation})). 

We thus wish to compute $\lim_{N \rightarrow
  \infty} \tilde{H}^{N,k,s,t}/N^{1+s+t}= \lim_{N\rightarrow \infty}
{H}^{N,k,s,t}/N^{1+s+t}$.

The proof of the equality (\ref{identityOne})-(\ref{identityTwo}) essentially follows from two
basic identities on formal series:
\begin{multline*}
\left(1-rX+r^2X^2-\cdots\right)^k\left(1-sX+s^2X^2-\cdots\right)\\
 = \sum_j (-1)^j \sum_{i=0}^j \binom{k+i-1}{i}r^is^{j-i}X^j
\end{multline*}
and
\begin{eqnarray}
\sum_{\substack{0\le j \le k-1 \\0\le q\le k-1}} (-1)^q \binom{k-1}{q}
q^j X^j &=& (-1)^{k+1} (k-1)!X^{k-1}. \label{cancellation}
\end{eqnarray}
We expand the definition of $\tilde{H}^{N,k,s,t}$ as a power series in
$q$. The first identity indicates that we should only look at the coefficient of
$q^{k-1}$, which we obtain by using the second identity (set $r :=
1/N$, $s:=1/(p+1)$). We then let $N$ tend to infinity, so the sum over
$p$ becomes a Riemann sum. Its limit is a $\beta$-integral, and thus a $\beta$-function appears, which can be expanded into a product of
$\Gamma$-functions, giving the first equality.

The equality (\ref{identityThree})-(\ref{identityFour}) is immediate
and is the only one to require the bound $k>t$.

For equality (\ref{identityTwo})-(\ref{identityThree})\footnote{This
  equality was first proved using \textsf{Mathematica}. Paul Abbott
  observed that the hypergeometric function that appears is
  Saalsch\"utzian and 
  extracted the following proof by tracing \textsf{Mathematica}'s output.}, we first
define 
\begin{eqnarray*}
H^{k,s,t}_a &=& k \sum_{i=0}^{\infty} \frac{\Gamma(k+a+i) 
  \Gamma(s+a+i+t+1)}{\Gamma(a+i+1)\Gamma(k+s+t+a+i+2)}\\
&=& \frac{k \Gamma (a+k) \Gamma (a+s+t+1)}{\Gamma (a+1) \Gamma (a+k+s+t+2)} \,_3F_2\left(\substack{1,a+k,a+s+t+1 \\ a+1,a+k+s+t+2};1\right),
\end{eqnarray*}
which satisfies $H^{k,s,t} = H^{k,s,t}_0 - H^{k,s,t}_{k-t}$. The
second equality is merely a consequence of the definition of
$\,_3F_2$. 

Since (see \cite{Mathematica})
$$
\,_3F_2 \left(\substack{1,c,d\\e,c+d-e+2};1\right) = \frac{c+d-e+1}{(c-e+1)(d-e+1)}\left(1-e+\frac{\Gamma(c+d-e+1)\Gamma(e)}{\Gamma(c)\Gamma(d)}\right),
$$
we get
$$
H^{k,s,t}_a = \frac{1}{1+s+t}\left(1-\frac{a\Gamma(a+k)\Gamma(a+s+t+1)}{\Gamma(a+1)\Gamma(a+k+s+t+1)}\right),
$$
which lets us prove equality (\ref{identityTwo})-(\ref{identityThree}) using the relation  $H^{k,s,t} = H^{k,s,t}_0 - H^{k,s,t}_{k-t}$.
\end{proof}

Let $G(\cdot)$ be the Barnes
$G$-function~\cite[Appendix]{HughesFirst}. It is a quick consequence of the Weyl dimension formula (see \cite[Equation~(18)]{BumpGamburd} that
$$
s_{\left<N^k\right>}\Ones{2k} \sim_N \fracG N^{k^2}.
$$

We use the previous Proposition to give a relatively concise expression
for $(\mathcal{M})(2k,r)$.
\begin{thm}
For a fixed $k \in \mathbb{N}$, the two series  
\begin{multline}
\sum_{r>0} (\mathcal{M})(2k,r) \frac{(\mathfrak{i}z)^r}{r!}
 \text{ \quad and \quad }
\fracG 
\sum_{d=1}^\infty
\sum_{\substack{\sort{s},\sort{t} \in \mathbb{N}^d}}\\
\left|\frac{ 1 }{ s_i!  t_j!(1+s_i+t_j)}\right|_{d\times d} 
\left| \frac{ H^{k,s_i,t_j} }{ s_i!  t_j!}
   \right|_{d\times d} z^{d+\total 
     (s_i+t_i)} \label{HAppear}
\end{multline}
\begin{multline}
\hfill
=
\fracG 
\sum_{d=1}^\infty
\sum_{\substack{\sort{s},\sort{t} \in \mathbb{N}^d}}
\left|\frac{ 1 }{ s_i!  t_j!(1+s_i+t_j)}\right|_{d\times d}^2 \hfill\\
\left(\prod_{i,j=1}^d
\frac{\Gamma(2k-t_j)\Gamma(k+s_i+1)}{\Gamma(k-t_j)\Gamma(2k+s_i+1)}\right)
z^{d+\total s_i+t_i}
\label{computational}
\end{multline}
\begin{multline} = \hspace{-0.6in}  \fracG 
\sum_{\substack{\lambda =\Frob{\sort{s}}{\sort{t}}\\\operatorname{rank} \lambda
  = d}}
s_\lambda \Ones{k}
\left| \frac{ \frac{\Gamma(2k-t_j)}{\Gamma(2k+s_i+1)} }{s_i! t_j!(1+s_i+t_j)}
   \right|_{d\times d}
   z^{|\lambda|}\label{SchurAppear}
\end{multline}
have equal coefficients of $z^r$ for $r < 2k+1$. For a fixed $r$, the
coefficients of $z^r$ for low values of $k$ can be meromorphically continued into each other.

Furthermore, by using Cauchy's Lemma, one can switch to an expression
involving products instead of determinants (i.e. a hypergeometric expression).
\label{concise}
\end{thm}

\begin{proof}
For Expression~(\ref{HAppear}), we proceed essentially by
substitution into Equation~(\ref{seriesDeterminant}), and looking at
terms of order $N^{k^2+r}$. Again, Cauchy's
Lemma is used repeatedly to reorganize determinants.

For Expressions~(\ref{computational}) or (\ref{SchurAppear}), we reorganized yet again the
determinants using Cauchy's Lemma into a form corresponding to
Formula~(\ref{SamraKingFormula}). We also summed over partitions
$\lambda$ instead of summing first over their rank $d$ then their
Frobenius coordinates $\sort{s},\sort{t}$. 

For a fixed $r$, both sides indeed admit meromorphic continuations in
$k$, which are equal by Carlson's Theorem\cite[Theorem 2.8.1,
p.~110]{AAR}\footnote{See also the author's
  thesis\cite{DehayeThesis}.}. Indeed, the LHS is shown to
admit a meromorphic continuation in $k$ using a Pochhammer
contour. The RHS' meromorphic continuation is already  written in
Expression~(\ref{computational}), if we admit that what is meant there
is the value of the meromorphic continuation in $k$ 
\emph{evaluated at $k$}. The difference of the two sides 
satisfies the hypotheses in Carlson's Theorem, in that its value is 0 at
integers, it is of exponential type, and type $<\pi$ along axes
parallel to the imaginary axis. The author has shown similar
statements in his thesis.

It is probably good to insist that the meromorphic continuation
of $\frac{\Gamma(2k-t_j)\Gamma(k+s_i+1)}{\Gamma(k-t_j)\Gamma(2k+s_i+1)}$
to the left has to be taken very carefully and cannot be obtained by
just plugging in values of $k$, once $k \le t$. We will discuss similar issues later, in Section~\ref{SectionConreyGhosh}.
\end{proof}

We now aim to replace the determinant left in Equation~(\ref{SchurAppear})
by a friendlier expression, a rational function of $k$.
\section{General shape of  $(\mathcal{M})(2k,r)$,
  $|\mathcal{M}|(2k,2h)$ and $|\mathcal{V}|(2k,2h)$} 
\label{together}
We now prove Theorem~\ref{thmComplexM}. 
\begin{proof}
By Equation~(\ref{HAppear}), we know that (for fixed $r$ and as
meromorphic functions of $k$) 
$$
\frac{\mathfrak{i}^r}{r!}(\mathcal{M})(2k,r) = \fracG \sum_{\substack{1 \le d \\ \sort{s},\sort{t}
  \in \mathbb{N}^d\\ d+\total (s_i+t_i) = r}} C(d,\sort{s},\sort{t})
\left| H^{k,s_i,t_j} \right|_{d \times d},  
$$
with $C(d,\sort{s},\sort{t}) \in \mathbb{Q}$, 
while for $s$ and $t$ fixed (and non-negative, of course),
Equation~(\ref{identityFour}) indicates that $H^{k,s,t}$ is a rational
function of $k$:
\begin{eqnarray}
H^{k,s,t} = \frac{1}{1+s+t} \prod_{i=-t}^s\frac{k+i}{2k+i}. \label{ratioH}
\end{eqnarray}
This already shows that we have a rational function of $k$ and that
the degree of its numerator equals the degree of its denominator. Equations~(\ref{ratioH}) and~(\ref{HAppear}) together, along
with the fact that $H^{k,s,t}=H^{-k,t,s}$ (a consequence of
Equation~(\ref{identityFour}),
explain why $X_r$ is even.

In order to determine the $Y_r$s a bit better, we need to investigate
possible denominators in the terms of $\left| H^{k,s_i,t_j} \right|_{d
  \times d}$. If $a$ is positive, 
$\left| H^{k,s_i,t_j} \right|_{d \times d}$ will have a factor of
$(2k+a)^{\alpha_a(r)}$ in its denominator if and only if $a$ is odd
 (because there is cancellation in
Formula~(\ref{ratioH})) and all of $s_1,\cdots, s_{\alpha_a}$ are greater than $a$. For this to happen, we need 
\begin{eqnarray}
r & = & d + \sum s_i + \sum t_i \notag\\
  & \ge & \alpha_a(r) + \sum_{i=1}^{\alpha_a(r)} (a+i-1) + \sum_{i=1}^{\alpha_a(r)} (i-1),\label{boundPole}
\end{eqnarray}
where the inequality is obtained by taking as small as possible values
for $d$
 (i.e. $\alpha_a(r)$), for the $s_i$'s (while requiring them to be
 different and greater or equal to $a$) and for the $t_i$'s (all different).

We turn this inequality around and get 
$$
\alpha_a(r) \le
\left\lfloor\frac{-a+\sqrt{a^2+4r}}{2}\right\rfloor.
$$

The case of $a$ negative is the same, exchanging the roles played
by $\sort{s}$ and~$\sort{t}$. 

Finally, the constant $D(r)$ ensuring that both $X_r$ and $Y_r$ are
monic can be found, thanks to Equations~(\ref{HAppear}) and
(\ref{ratioH}), taking $\lim_{k \rightarrow \infty}$: 
\begin{eqnarray}
D(r) =  \sum_{\substack{1 \le d \\ \sort{s},\sort{t}
  \in \mathbb{N}^d\\ d+\total (s_i+t_i) = r}}
\left|\frac{1}{s_i!t_j!(1+s_i+t_j)}\right|^2_{d \times d}
\frac{1}{2^{d+\total (s_i+t_i)}}=\frac{1}{r!2^r},
\label{BorodinTag}
\end{eqnarray}
where this last equality is left to the reader. 

Actually, this last equality is enough to also guarantee that $X_r(u)$ and
$Y_r(u)$ both have integer coefficients: just substitute for $H^{k,s,t}$ in
Equation~(\ref{HAppear}) 
$$
H^{k,s,t} = \frac{1}{1+s+t} \prod_{i=-t}^s\frac{1}{2k+i}
\left(k\sum_{i=0}^{s+t} h_i k^{i}\right)
$$
for the appropriate (integer) $h_i$s (in particular, $h_{s+t}=1$).

This proves Equation~(\ref{polComplex}), at least for large $k$.

Meromorphic continuation has already been obtained in Theorem~\ref{concise}.
\end{proof}

\begin{thm}
\label{otherMoments}
For $h\in \mathbb{N}$, there are polynomials $\tilde{X}_{2h},\hat{X}_{2h}, $ with
  integer coefficients and $\deg \hat{X}_{2h} = \deg X_{2h} > \deg
  \tilde{X}_{2h}$ such that as meromorphic functions of $k$, 
\begin{eqnarray*}
|\mathcal{M}|(2k,2h) &=& \hat{C}(h) \fracG
\frac{\hat{X}_{2h}(2k)}{Y_{2h}(2k)},\\
|\mathcal{V}|(2k,2h) &=& \tilde{C}(h) \fracG
\frac{\tilde{X}_{2h}(2k)}{Y_{2h}(2k)},
\end{eqnarray*}
where $Y_{r}(u)$ is as defined in Theorem~\ref{thmComplexM}.

Moreover (but this is conjectural), the numerators are additionally 
monic polynomials\footnote{This is
the normalization we will keep later, when discussing data about those
polynomials.} when $\hat{C}(h) = \frac{1}{2^{2h}}$ and
$\tilde{C}(h)=\frac{(2h)!}{h! 2^{3h}}$, and $ \deg X_{2h} - \deg
  \tilde{X}_{2h}= 2h$.
\end{thm}
\begin{proof}
For fixed integer $r$ and large integer $k$, most of this follows immediately from Equations~(\ref{relationTwo})
and (\ref{relationFour}), combined with Theorem~\ref{thmComplexM}. 

The fact that $\deg \tilde{X}_{2h}<\deg X_{2h}$ for instance is a consequence of 
$$
\left(\mathcal{M}\right)(2k,r) \sim_k
\left(-\frac{\mathfrak{i}}{2}\right)^r \fracG,
$$
which we use in Equation~(\ref{eqnNotContinuation}):
$$
\sum_{j=0}^{2h} \binom{2h}{j}
\left(\frac{\mathfrak{i}}{2}\right)^{j} \left(-\frac{\mathfrak{i}}{2}\right)^{2h-j}=0.
$$

We can similarly show that \emph{if it exists}, $\hat{C}(h) =
\frac{1}{2^{2h}}$. The constant $\tilde{C}(h)$ is more mysterious, and
involves lower order terms in $k$ of Equation~(\ref{ratioH}).

The meromorphic continuation is obtained as in the proof of Theorem~\ref{thmComplexM}.
\end{proof}

\textbf{Remark. } Unfortunately, within their degree restrictions, the
$X_{r}(u), \tilde{X}_{2h}(u)$ and $\hat{X}_{2h}(u)$ polynomials still
look utterly random. We merely have an expression for them as a sum of determinants of rank $d\le\sqrt{r}$ (resp.~$2h$). This expression is relatively quick and
allows at least to compute a few of those polynomials. 

\section{Computational data}
\label{data}
\subsection{The polynomials $X_r(u)$, $\tilde{X}_{2h}(u)$ and $\hat{X}_{2h}(u)$}
We present our data for $(\mathcal{M})(2k,r)$ first, in
Table~\ref{dataComplexM}, followed by data on $|\mathcal{M}|(2k,2h)$
in Table~\ref{dataNormM} and finally on $|\mathcal{V}|(2k,2h)$ in
Table~\ref{dataNormV}. Everything extends numerical results
previously published, for instance in \cite{Hall2004,Hall2002} (but
those rely on \cite{HughesJoint}) or \cite{CRS} (which is limited to
$k=h$).
\ForArxiv{Extended versions of those tables are also made available in
  the source of this \arXiv\, submission or (possibly more) at \cite{DehayeData}.}
\ForJournal{Extended versions of those tables are also made available in
 the expanded content or (possibly more) at \cite{DehayeData}.}

To obtain those tables, we have implemented
Equation~(\ref{computational}), which is the most computationally
accessible version of the formulas available in
Theorem~\ref{concise}. \ForArxiv{A \textsf{Magma} implementation of
  this algorithm is also part of this \arXiv\, submission.}
\ForJournal{A \textsf{Magma} implementation of
  this algorithm is also available as expanded content.}

\begin{sidewaystable}
\caption{
The first polynomials $X_r(u)$, i.e the numerators in
  $(\mathcal{M})(u,r)$. \ForJournal{Data up to $r=60$ available as
    expanded content or in~\cite{DehayeData}.} \ForArxiv{Data up to $r=60$ available attached to the source
    of this \arXiv\,submission or at \cite{DehayeData}.}  
}
\label{dataComplexM}
\begin{tabular}{c|p{0.9\linewidth}}
$r$& $X_r(u)$\\ \hline
1&$
1
$\\\hline
2&$
u^{2}
-2
$\\\hline
3&$
u^{2}
-4
$\\\hline
4&$
u^{4}
-16u^{2}
+66
$\\\hline
5&$
u^{4}
-20u^{2}
+114
$\\\hline
6&$
u^{8}
-51u^{6}
+864u^{4}
-5554u^{2}
+4860
$\\\hline
7&$
u^{8}
-57u^{6}
+1134u^{4}
-8758u^{2}
+8520
$\\\hline
8&$
u^{10}
-113u^{8}
+4620u^{6}
-86332u^{4}
+682844u^{2}
-765660
$\\\hline
9&$
u^{10}
-121u^{8}
+5460u^{6}
-115564u^{4}
+1053964u^{2}
-1457820
$\\\hline
10&$
u^{14}
-220u^{12}
+18897u^{10}
-831010u^{8}
+20196928u^{6}
-260164440u^{4}
+1428629724u^{2}
-2060092440
$\\\hline
11&$
u^{14}
-230u^{12}
+20997u^{10}
-996820u^{8}
+26447168u^{6}
-374214600u^{4}
+2270621484u^{2}
-3994446960
$\\\hline
12&$
u^{18}
-363u^{16}
+52929u^{14}
-4083011u^{12}
+183649422u^{10}
-4906031274u^{8}
+73323636100u^{6}
-512994314412u^{4}
+1371835414728u^{2}
-927651213720
$\\\hline
13&$
u^{18}
-375u^{16}
+57141u^{14}
-4663655u^{12}
+224398746u^{10}
-6467410170u^{8}
+105010072036u^{6}
-806857605660u^{4}
+2461218471576u^{2}
-1755890884440
$\\\hline
14&$
u^{22}
-582u^{20}
+141344u^{18}
-18977780u^{16}
+1571817537u^{14}
-84339778978u^{12}
+2962887441370u^{10}
-66386724069396u^{8}
+884603961264548u^{6}
-6212383525692744u^{4}
+19176051246319080u^{2}
-13863690471430800
$\\\hline
15&$
u^{22}
-596u^{20}
+149296u^{18}
-20838716u^{16}
+1807941481u^{14}
-102286957136u^{12}
+3809004157906u^{10}
-90891702433976u^{8}
+1298188100828836u^{6}
-9917808021410976u^{4}
+33986748108863880u^{2}
-25682708695644000
$\\\hline
16&$
u^{24}
-836u^{22}
+295486u^{20}
-58491716u^{18}
+7245863641u^{16}
-593291868896u^{14}
+32861804018536u^{12}
-1227084273320096u^{10}
+29900504376591736u^{8}
-444180655702337856u^{6}
+3616035044845449600u^{4}
-13500165816324763200u^{2}
+10671545982659562000
$\\\hline
17&$
u^{24}
-852u^{22}
+308606u^{20}
-62999492u^{18}
+8101703961u^{16}
-692989945072u^{14}
+40321523165416u^{12}
-1589469869122752u^{10}
+41098203910503416u^{8}
-652694167393180032u^{6}
+5757854141711318400u^{4}
-23590053001525406400u^{2}
+19761261673907754000
$\\\hline
18&$
u^{30}
-1216u^{28}
+641547u^{26}
-195081042u^{24}
+38335269063u^{22}
-5171814422892u^{20}
+495753742037253u^{18}
-34353739684203042u^{16}
+1726507702228490928u^{14}
-62290017635596811632u^{12}
+1575250938092261972152u^{10}
-26886933063310680515376u^{8}
+293595553738705518511056u^{6}
-1882598606626601433513600u^{4}
+5855125431247144869877200u^{2}
-4699357338080820827412000
$\\\hline
19&$
u^{30}
-1234u^{28}
+663111u^{26}
-206226048u^{24}
+41629109007u^{22}
-5794171874298u^{20}
+575320671855777u^{18}
-41443936954862628u^{16}
+2171988993390059952u^{14}
-81956498940701768368u^{12}
+2174685878160406187416u^{10}
-39111313358222167862304u^{8}
+452360970074645727302736u^{6}
-3084756281794829726025120u^{4}
+10210913321050424344698000u^{2}
-8861284072193198189544000
$\\\hline
20&$
u^{34}
-1615u^{32}
+1140143u^{30}
-467224385u^{28}
+124593557421u^{26}
-22981261798995u^{24}
+3040237566735165u^{22}
-294611133821587635u^{20}
+21107532245623967310u^{18}
-1116405478738744697410u^{16}
+43058312795636550000904u^{14}
-1183070247664529791035320u^{12}
+22374172979188549647921632u^{10}
-277662183945403036368852000u^{8}
+2095071747708073688848702224u^{6}
-8269151494407104768839910640u^{4}
+12529695816553717113566335200u^{2}
-6099189940914050054558484000
$

\end{tabular}
\end{sidewaystable}

\begin{sidewaystable}[p]
\caption{
The first polynomials $\hat{X}_{2h}(u)$, i.e. the numerators
  in $|\mathcal{M}|(u,2h)$.
\ForJournal{Data up to $h=30$ available as
    expanded content or in~\cite{DehayeData}.}
\ForArxiv{Data up to $h=30$ available attached to the source of this
\arXiv\,submission or at \cite{DehayeData}.}
}  \label{dataNormM}
\begin{tabular}{c|p{0.9\linewidth}}
$r$& $\hat{X}_r(u)$\\ \hline
2&$
u^{2}
$\\\hline
4&$
u^{4}
-8u^{2}
-6
$\\\hline
6&$
u^{8}
-33u^{6}
+198u^{4}
+74u^{2}
-360
$\\\hline
8&$
u^{10}
-81u^{8}
+1740u^{6}
-8284u^{4}
-7716u^{2}
+34020
$\\\hline
10&$
u^{14}
-170u^{12}
+9597u^{10}
-215560u^{8}
+1846928u^{6}
-4247400u^{4}
-12317076u^{2}
+42366240
$\\\hline
12&$
u^{18}
-291u^{16}
+30177u^{14}
-1379507u^{12}
+28177518u^{10}
-236602818u^{8}
+604630084u^{6}
+1570591476u^{4}
-10008266040u^{2}
+7829929800
$\\\hline
14&$
u^{22}
-484u^{20}
+90384u^{18}
-8378492u^{16}
+415889897u^{14}
-11196067680u^{12}
+157699171570u^{10}
-1023611526808u^{8}
+1699483809828u^{6}
+11589901952544u^{4}
-62361799232760u^{2}
+44754182272800
$\\\hline
16&$
u^{24}
-708u^{22}
+198590u^{20}
-28525892u^{18}
+2275085529u^{16}
-102837376096u^{14}
+2598141390568u^{12}
-34807690054560u^{10}
+213458763180152u^{8}
-261862022455104u^{6}
-3402805264433280u^{4}
+19256263380043200u^{2}
-11718802173078000
$\\\hline
18&$
u^{30}
-1054u^{28}
+460431u^{26}
-109299828u^{24}
+15577804767u^{22}
-1394331670638u^{20}
+79872695247657u^{18}
-2932723486507728u^{16}
+68022586503825552u^{14}
-962308385613255088u^{12}
+7682283932820069016u^{10}
-26475220331016986304u^{8}
-59889950570120914224u^{6}
+976582356673028315040u^{4}
-3441287004848413282800u^{2}
+1366282646437284576000
$\\\hline
20&$
u^{34}
-1415u^{32}
+840943u^{30}
-275540385u^{28}
+55049482221u^{26}
-7022476724835u^{24}
+584090828573565u^{22}
-31869278744265555u^{20}
+1134427249824868110u^{18}
-25880772100948222330u^{16}
+365485578445889268104u^{14}
-2970099871666499086840u^{12}
+10773785732163438366432u^{10}
+24904735536575464181280u^{8}
-474478390713139651278576u^{6}
+1993984711160163968152080u^{4}
-1770512318771949573760800u^{2}
+214967318998766249916000
$

\end{tabular}
\end{sidewaystable}

\begin{sidewaystable}[p]
\caption{
The first polynomials $\tilde{X}_{2h}(u)$, i.e. the numerators
  in $|\mathcal{V}|(u,2h)$. \ForJournal{Data up to $h=30$ available  as
    expanded content or in~\cite{DehayeData}.}\ForArxiv{Data up to $h=30$ available attached to the source
    of this \arXiv\,submission or at \cite{DehayeData}.} } \label{dataNormV}
\begin{tabular}{c|p{0.9\linewidth}}
$r$& $\tilde{X}_r(u)$\\ \hline
2&
$
1
$\\\hline
4&
$
1
$\\\hline
6&
$
u^{2}
-9
$\\\hline
8&
$
u^{2}
-33
$\\\hline
10&
$
u^{4}
-90u^{2}
+1497
$\\\hline
12&
$
u^{6}
-171u^{4}
+6867u^{2}
-27177
$\\\hline
14&
$
u^{8}
-316u^{6}
+30702u^{4}
-982572u^{2}
+6973305
$\\\hline
16&
$
u^{8}
-484u^{6}
+76902u^{4}
-4461348u^{2}
+67692705
$\\\hline
18&
$
u^{12}
-766u^{10}
+215847u^{8}
-27766980u^{6}
+1653656895u^{4}
-41530140126u^{2}
+337968054585
$\\\hline
20&
$
u^{14}
-1055u^{12}
+421093u^{10}
-79486155u^{8}
+7242179715u^{6}
-290444510205u^{4}
+4099101803991u^{2}
-8381907513945
$\\\hline
22&
$
u^{16}
-1496u^{14}
+892108u^{12}
-272180808u^{10}
+45430344630u^{8}
-4121412379560u^{6}
+189676636728876u^{4}
-3674923533427896u^{2}
+14539253947899345
$\\\hline
24&
$
u^{18}
-1961u^{16}
+1566628u^{14}
-658984788u^{12}
+157743552510u^{10}
-21750520014270u^{8}
+1678578114026196u^{6}
-67707100461703716u^{4}
+1235110338400818825u^{2}
-6787336148294472225
$\\\hline
26&
$
u^{20}
-2610u^{18}
+2860437u^{16}
-1718473240u^{14}
+620475009522u^{12}
-139083336332460u^{10}
+19348398203611266u^{8}
-1624490941247619480u^{6}
+77190294570345945549u^{4}
-1813095317449668401010u^{2}
+15009483262024846096425
$\\\hline
28&
$
u^{22}
-3243u^{20}
+4462647u^{18}
-3407674501u^{16}
+1586340567882u^{14}
-466277764083726u^{12}
+86845227411024846u^{10}
-10042821279688179978u^{8}
+688582088681764130469u^{6}
-25698037955845496067927u^{4}
+444470604942195922015755u^{2}
-2654155080367803900605025
$\\\hline
30&
$
u^{28}
-4190u^{26}
+7631083u^{24}
-7953124300u^{22}
+5258554468937u^{20}
-2313326757869890u^{18}
+691451285514065259u^{16}
-141062107217586416040u^{14}
+19477099336547993586171u^{12}
-1781103872658227723795970u^{10}
+103764470143371018680338137u^{8}
-3607131084573924222894990540u^{6}
+66647887693999747894954784187u^{4}
-515514421669410774166185623070u^{2}
+658183121944091618062137174225
$

\end{tabular}
\end{sidewaystable}

\subsection{The roots of $\tilde{X}_{2h}(u)$}
\label{position}
It has been suggested before, based on limited numerical data, that
the polynomials $\tilde{X}_{2h}(u)$ have only real roots. In fact we list in
Table~\ref{rootcount} the number of real roots and degree for each
such polynomial.
\begin{table}
\caption{\label{rootcount} The degree and the number of real roots of
  $\tilde{X}_{2h}$. The $h$s for which there are non-real roots are highlighted.}
\begin{tabular}{|c||c|c|c|c|c|c|c|c|c|c|c}
\hline
h&1&2&3&4&5&6&7&8&9&10&
\\
\hline
\hline
$\deg(\tilde{X}_{2h})$&
0& 0& 2& 2& 4& 6& 8& 8& 12&14&
\\
\hline
$\#$ real roots&
0& 0& 2& 2& 4& 6& 8& 8& 12&14&
\\
\hline
\multicolumn{11}{c}{}\\
\hline
h&11&12&13&14&15&16&17&18&19&20&
\\
\hline
\hline
$\deg(\tilde{X}_{2h})$&
 16& 18& 20& 22& 28& 28& 30& 34& 36& 38&
\\
\hline
$\#$ real roots&
 16& 18& 20& 22& 28& 28& 30& 34& 36& 38&
\\
\hline
\multicolumn{11}{c}{}\\
\hline
h&\textbf{21}&22&\textbf{23}&\textbf{24}&25&\textbf{26}&\textbf{27}&\textbf{28}&\textbf{29}&\textbf{30}&
\\
\hline
\hline
$\deg(\tilde{X}_{2h})$& \textbf{44}& 46& \textbf{48}& \textbf{50}& 54& \textbf{56}& \textbf{62}& \textbf{64}& \textbf{66}& \textbf{72}&\\
\hline
$\#$ real roots&\textbf{40}& 46& \textbf{44}& \textbf{46}& 54&
\textbf{52}& \textbf{58}& \textbf{60}& \textbf{62}& \textbf{68}&\\
\hline
\end{tabular}
\end{table}
One quickly observes that $\tilde{X}_{42}(u)$ (of course!) is actually
the first polynomial to break the initial fluke and have non-real
roots. It is not clear at this point if this is related to a similar 
observation on the last line of \cite{Hall2002} and throughout \cite{Hall2004}.

{
\tiny
\begin{align*}
X_{42}(u)=u^{44} &- 12302u^{42} + 69239935u^{40} - 236610412148u^{38}\\ &+ 549459541784707u^{36} - 919748248913270486u^{34} + 1148989069656897835213u^{32}\\& - 1094474723973849448826480u^{30} + 
    805533314533281755701371226u^{28} \\&-
    461541928967718110253944237052u^{26} +
    206514429127544387915748094513446u^{24} \\&-
    72119441118339869972121541587076920u^{22} \\& + 
    19577196457693502603026719624834404502u^{20}\\& - 4099121776759328236737053383626986012604u^{18}\\& + 654170727960937096861203148250462720819850u^{16}\\& - 
    78212503734767115379758317319774926243800176u^{14} \\&+ 6836980008003428572296900814856434321006155189u^{12}\\& - 422028250886223501142365592098345343850710857462u^{10}\\ &+ 
    17476800084974190439148752639441918166326024419531u^8 \\&- 448540393629268182677088044978029477583305447285620u^6 \\&+ 6253526937210642323596984565394593401672539709730775u^4 \\&- 
    37013087756228993438266827460643377762894550851248750u^2 \\&+
    36216052456609571501642100973941635690472733838765625.
\end{align*}
}
This polynomial has four non-real roots  ($\pm 18.8631835 \pm
0.0090603\mathfrak{i}$) that show up at once, since they would have to come in
pairs of conjugate pairs by evenness of $\tilde{X}_{2h}(u)$. One could
wonder why non-real roots show up so late, and if there is actually a
good reason for this.  

\begin{fact}
The polynomials $\tilde{X}_{2h}(u)$ tend to have many, \emph{but not all}, of
their roots real. For instance, for high $h$, $\tilde{X}_{2h}(u)$ has one root
very close by to every odd integer between $h$ and $2h$.
\end{fact}

\begin{sidewaysfigure}
\includegraphics[width=\textwidth]{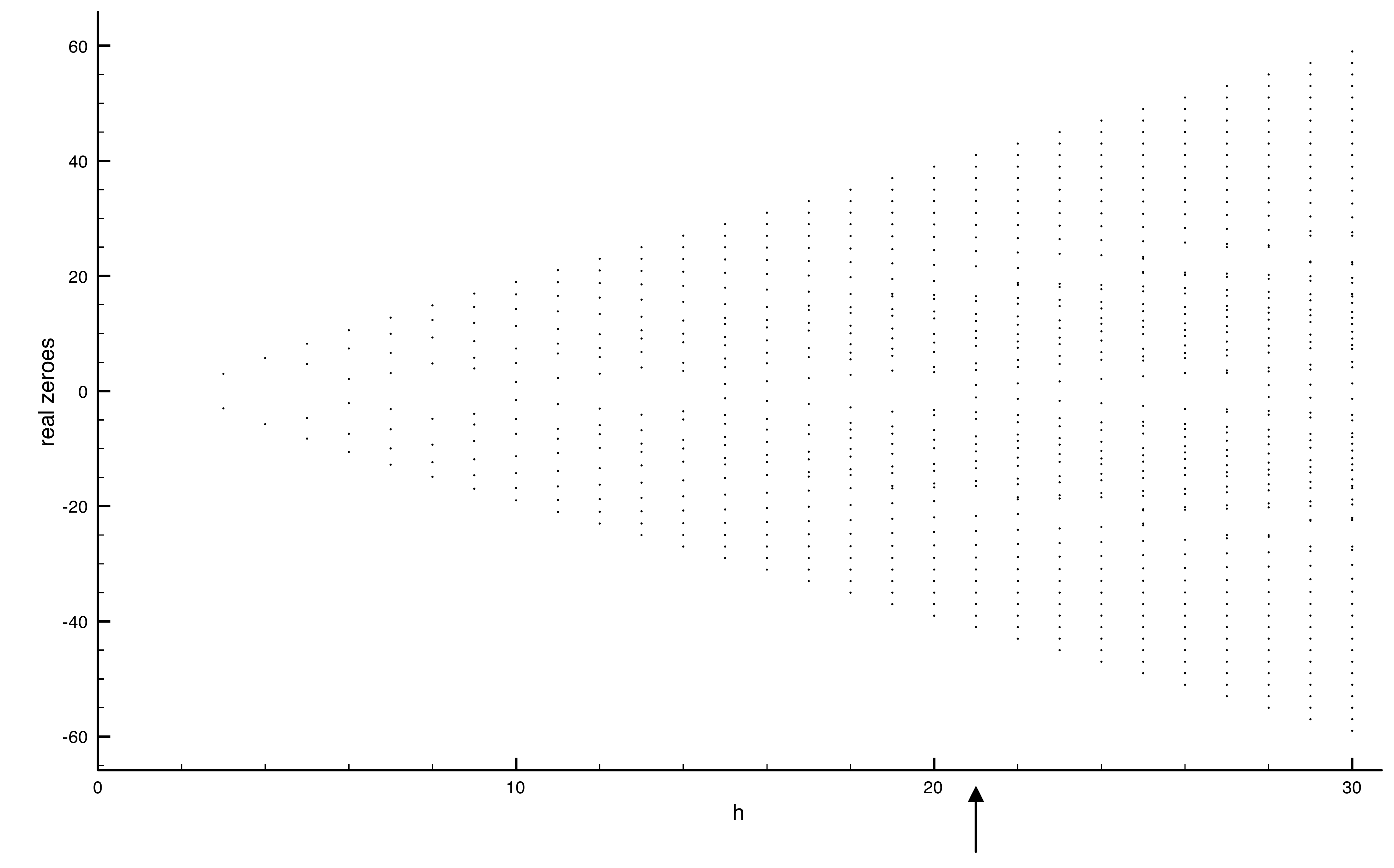}
\caption{
The roots of $\tilde{X}_{2h}(u)$. The line corresponding to
  $\tilde{X}_{42}(u)$, i.e. where the first real zeroes go missing, has been indicated.}
\label{fig_one}
\end{sidewaysfigure}

We first present graphical clues for this fact in Figure~\ref{fig_one}, which
depicts the position of the real roots for $h=1$ to $h=30$. It thus
omits the complex roots.


We now explain the fact. It helps at this point to remember that
$\tilde{X}_{2h}(u)$ is obtained by summing various $X_r(u)$ for $r\le
2h$, which are themselves obtained from Equation~(\ref{HAppear}), for
instance. Furthermore, the summand in that Equation associated to
$d,\vect{s},\vect{t}$ (with $r=d+\sum_i s_i +t_i$) will have poles (as
a function of $u=2k$) at the odd integers $a$ such that $-s_1 \le a
\le t_1$ (this uses Lemma~\ref{CauchysLemma} to expand the determinant
in $H^{k,s_i,t_j}$s). For each pole $a$, there are a few summand where
this pole 
comes with multiplicity exactly $\alpha_a(r)$, but for most others the
multiplicity is lower (see Equation~(\ref{boundPole})). So if we sum all of those terms, and multiply
by $Y_{2h}(u)$ (the common denominator) to obtain $\tilde{X}_{2h}(u)$,
a vast majority of terms factor a $(u-a)$ out. We thus have an
expression of the form
$$
\tilde{X}_{2h}(u) = (u-a) P_1(u) + P_2(u),
$$
where the coefficients of $P_1(u)$ are expected to be much bigger than
the coefficients of $P_2(u)$ (simply because much more terms are summed to
obtain $P_1(u)$ than $P_2(u)$). Hence, we should expect
$\tilde{X}_{2h}(u)$ to change sign when $u$ travels along the real
axis from below $a$ to above $a$ (because $|P_1(a)|>\!\!|P_2(a)|$ and
$(u-a)$ changes sign)  and we
know that a root will be around $u=a$. This is especially true if
$a>r/2$, because the restrictions impose then $s_1>a>s_2$, and as a consequence
$\alpha_r(a)=1$ and the phenomenon described just now is
accentuated. We present in Table~\ref{largestRoot} some numerical data
associated to this phenomenon.

It is obvious from Figure~\ref{fig_one} that a lot is yet to be understood about
the polynomials $\tilde{X}_{2h}(u)$. For instance, it is not clear if
asymptotically in $h$ there is a positive proportion of real roots. 
\begin{table}
\caption{The largest root of $\tilde{X}_{2h}(u)$}
\begin{tabular}{r|r|r|r}
h& largest root of $\tilde{X}_{2h}(u)$ & difference with $2h-1$& log. difference \\
\hline
1& no root & no root&no root\\
2& no root & no root&no root\\
3  &  3.0000000000000000000  &  2.000000000  &  0.69315 \\
4  &  5.7445626465380286598  &  1.255437354  &  0.22748 \\
5  &  8.2448923938491831987  &  0.7551076062  &  -0.28090 \\
6  &  10.568920444013080343  &  0.4310795560  &  -0.84146 \\
7  &  12.769459455674733521  &  0.2305405443  &  -1.4673 \\
8  &  14.886048429155973920  &  0.1139515708  &  -2.1720 \\
9  &  16.948550444560344620  &  0.05144955544  &  -2.9672 \\
10  &  18.978943770872905688  &  0.02105622913  &  -3.8606 \\
11  &  20.992206162055831068  &  0.007793837944  &  -4.8544 \\
12  &  22.997383184072186530  &  0.002616815928  &  -5.9458 \\
13  &  24.999198051064882757  &  0.0008019489351  &  -7.1285 \\
14  &  26.999774030173017860  &  0.0002259698270  &  -8.3951 \\
15  &  28.999941044846106152  & $ 5.895515389 \times 10^{-5}$  &  -9.7388 \\
16  &  30.999985671005722891  & $ 1.432899428\times 10^{-5}$  &  -11.153 \\
17  &  32.999996738730003824  & $ 3.261269996\times 10^{-6}$  &  -12.633 \\
18  &  34.999999301847217917  & $ 6.981527821\times 10^{-7}$  &  -14.175 \\
19  &  36.999999858891343014  & $ 1.411086570\times 10^{-7}$  &  -15.774 \\
20  &  38.999999972983353984  & $ 2.701664602\times 10^{-8}$  &  -17.427 \\
21  &  40.999999995085836086  & $ 4.914163914\times 10^{-9}$  &  -19.131 \\
22  &  42.999999999148595422  & $ 8.514045781\times 10^{-10}$  &  -20.884 \\
23  &  44.999999999859167358  & $ 1.408326421\times 10^{-10}$  &  -22.683 \\
24  &  46.999999999977712180  & $ 2.228782021\times 10^{-11}$  &  -24.527 \\
25  &  48.999999999996618870  & $ 3.381129731\times 10^{-12}$  &  -26.413 \\
26  &  50.999999999999507453  & $ 4.925468142\times 10^{-13}$  &  -28.339 \\
27  &  52.999999999999930988  & $ 6.901186254\times 10^{-14}$  &  -30.304 \\
28  &  54.999999999999990686  & $ 9.313971788\times 10^{-15}$  &  -32.307 \\
29  &  56.999999999999998787  & $ 1.212486889\times 10^{-15}$  &  -34.346 \\
30  &  58.999999999999999847  & $ 1.524414999\times 10^{-16}$  &  -36.420 \\
\end{tabular}
\label{largestRoot}
\end{table}
\section{Alternative expressions}
\label{alternative}
\subsection{Using Macdonald's ninth variation of Schur functions}
Define, as in \cite{Noumi1} and \cite{NoumiBook}, and similarly to \cite{MacdonaldNinth},
\begin{eqnarray}
\label{variation}
\tilde{s}_\lambda^{(R)}:= \left| \tilde{h}^{(R-j+1)}_{\lambda_i-i+j} \right|_{l(\lambda) \times l (\lambda)},
\end{eqnarray}
with
$$
\tilde{h}^{(R)}_{k} := \frac{(R-1)!}{(R+k-1)! k!}.
$$
We first prove that this variation of Schur functions satisfies a \emph{Giambelli identity}.
\begin{prop}
Let $\lambda$ be a partition and $\Frob{\sort{s}}{\sort{t}}$ its Frobenius coordinates, of rank~$d$. Then,
$$
\tilde{s}_\lambda^{(R)}= \left| \tilde{s}^{(R)}_{(s_i|t_j)} \right|_{d
  \times d} =
\left| \frac{ \frac{\Gamma(R-t_j)}{\Gamma(R+s_i+1)} }{s_i! t_j!(1+s_i+t_j)}
   \right|_{d\times d}.
$$
\end{prop}
Note how this provides a second determinantal expression for this
variation of Schur functions, but with a matrix of different rank.
\begin{proof}
We intend to use Exercise 3.21 in Macdonald's book, but to show that the exercise applies, we need to prove:
\begin{eqnarray*}
\tilde{s}^{(R)}_{(p|q)} &:= &\det \begin{pmatrix}
\tilde{h}^{(R)}_{p+1} & \tilde{h}^{(R-1)}_{p+2} & \cdots & 
\cdots & \cdots & \tilde{h}^{(R-q)}_{p+q+1}\\
1 & \tilde{h}^{(R-1)}_{1} & \tilde{h}^{(R-2)}_{2} & \cdots & \cdots  &
\tilde{h}^{(R-q)}_{q}\\
0 & 1 & \tilde{h}^{(R-2)}_{1} & \tilde{h}^{(R-3)}_{2} & \cdots & 
\tilde{h}^{(R-q)}_{q-1}\\
\vdots & \ddots & \ddots & \ddots & \ddots&  \vdots\\
\vdots & \ddots & \ddots & \ddots & \ddots&  \vdots\\
0 & \cdots & \cdots &   0 & 1 & \tilde{h}^{(R-q)}_{1}\\
\end{pmatrix}_{(q+1)\times (q+1)}\\
& =& \frac{1}{p!q!(1+p+q)}\frac{\Gamma(R-q)}{\Gamma(R+p+1)}. 
\end{eqnarray*}
This can be shown by expanding the determinant along the last column
to obtain 
$$ \tilde{s}^{(R)}_{(p|q)} = (-1)^q \tilde{h}^{(R-q)}_{p+q+1} +
\sum_{i=1}^q (-1)^{i+1}\tilde{h}^{(R-q)}_i
\tilde{s}^{(R)}_{(p|q-i)}.$$
Subtract the LHS from the RHS, proceed by induction on $q$, factor
out $\frac{\Gamma(R-q)}{\Gamma(R+p+1)}$ and the
result then follows from the following equalities, for $p$ and $q$
positive integers:
\begin{multline*}
 \frac{(-1)^q}{(p+q+1)!} - \sum_{i=1}^q
\frac{(-1)^{i}}{i!p!(q-i)!(p+q-i+1)} - \frac{1}{p!q!(1+p+q)} \\
=\frac{(-1)^q}{(p+q+1)!} + \frac{p+q+q \quad _2F_1(\substack{1-q\quad -p-q\\1-p-q};1)}{(p+q)(p+q+1)p!q!}- \frac{1}{p!q!(1+p+q)}\\
= \frac{(-1)^q}{(p+q+1)!} + \frac{q \quad _2F_1(\substack{1-q\quad
    -p-q\\1-p-q};1)}{(p+q)(p+q+1)p!q!}\\
=\frac{(-1)^q}{(p+q+1)!}  
+\frac{q}{(p+q)(p+q+1)p!q!}\frac{(q-1)!}{(1-p-q)^{(q+1)}}=0,
\end{multline*}
the last one being a consequence of Gauss's Hypergeometric theorem.

The theorem now results directly from Exercise 3.21 in Macdonald's book.
\end{proof}
In essence, this Proposition allows us to switch from a Giambelli-type
expression to a Jacobi-Trudi expression. It immediately leads to a
simplified version of Theorem~\ref{concise}.
\begin{thm}
\label{ninthTheorem}
With $G(\cdot)$ the Barnes $G$-function, and
  $\tilde{s}_\lambda$ defined as in Equation~(\ref{variation}), 
\begin{eqnarray*}
\sum_{r>0} (\mathcal{M})(2k,r) \frac{(\mathfrak{i}z)^r}{r!}
&=& \fracG
\sum_\lambda
\tilde{s}^{(2k)}_\lambda
s_\lambda \Ones{k}
   z^{|\lambda|},
\end{eqnarray*}
in the sense that their coefficients of $z^r$  are equal for fixed $r$ and large
enough $k$ so the coefficient in the LHS is defined. 
\end{thm}
\subsection{Imitating the Cauchy identity}
\label{Cauchy}
We can also give an alternative for the expression in
Equation~(\ref{SchurAppear}), proceeding as in
Gessel's theorem in its lead up to the Cauchy
identity (see~\cite{TWlongest}). This uses Theorem~\ref{ninthTheorem}.
\begin{thm} 
\begin{multline*}
\fracG \sum_{r>0} (\mathcal{M})(2k,r)
\frac{(\mathfrak{i}z)^r}{r!} \\
\begin{split}
& =&&  \lim_{n \rightarrow \infty} \det \left( 
\left( h_{j-i}\Ones{k} \right)_{n \times \infty}
\cdot
\left( \tilde{h}_{i-j}^{(2k-n+j)} z^{i-j}\right)_{\infty \times n}
\right)\\
& =&&   \lim_{n \rightarrow \infty} \det \left( 
\left( h_{j-i}\Ones{k} z^{j-i}\right)_{n \times \infty}
\cdot
\left( \tilde{h}_{i-j}^{(2k-n+j)}\right)_{\infty \times n}
\right)\\
&=&& \lim_{n \rightarrow \infty} \det \left( \sum_{l \ge 0} h_{l-i}
  \Ones{k} \tilde{h}_{l-j}^{(2k-n+j)} z^{l-j} \right)_{n \times n}\\
&=&&  \lim_{n \rightarrow \infty} \det \left( \sum_{l \ge 0}
  \binom{l-i+k-1}{k-1} \frac{(2k-n+j-1)!}{(l-j)!(2k-n+l-1)!} z^{l-j}
\right)_{n \times n}
,
\end{split}
\end{multline*}
in the sense that their coefficients of $z^r$ are equal for fixed $r$
and large enough $k$ so the coefficient in the LHS is defined.  
The factorials on the last line should really be evaluated in groups, to
give 0 if $l < j$, and $\frac{\Gamma(2k-n+j)}{\Gamma(2k-n+l)(l-j)!}$ otherwise.
\end{thm}

Note that this can be truncated significantly when we are
after only $$\sum_{0<r\le S} (\mathcal{M})(2k,r)
\frac{(\mathfrak{i}z)^r}{r!}$$ for a finite $S$
(i.e. when we are computing the head of the sequence of polynomials):
we can drop the limit in $n$ and settle for a sufficiently big $n$ instead, and
then cut the matrices in their infinite directions as well. 

In Gessel's Theorem, in order to get to the other side of the Cauchy
identity, one would then observe that the matrix on the last line is
Toeplitz, and then use Szeg\"o's theorem. Of course, that fails here
because the matrix on the last line is not Toeplitz.

\section{The result of Conrey and Ghosh}
\label{SectionConreyGhosh}
As explained in the introduction, Conrey and Ghosh's theorem \cite{ConreyGhosh} that
$\mathcal{J}(2,1) = \frac{e^2-5}{4 \pi}$ immediately leads to a
conjecture that $|\mathcal{V}|(2,1) = \frac{e^2-5}{4 \pi} $ as
well. Our main concern is that we only know $|\mathcal{V}|(2k,2h)$,
for integer $h$, through Equations~(\ref{relationTwo}) and
(\ref{relationFour}) (while we would need $h=1/2$). 

\begin{figure}
\includegraphics[width=\textwidth]{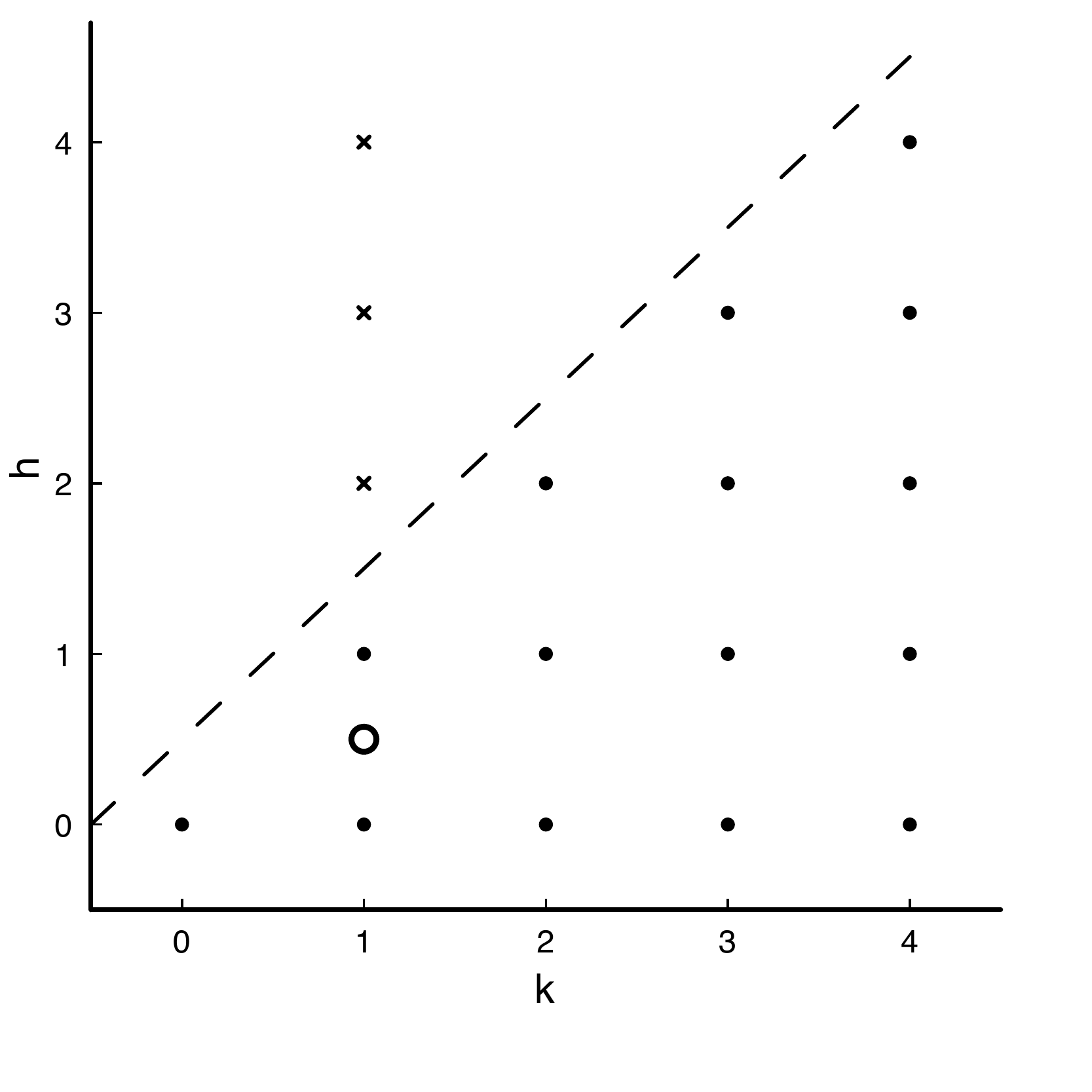}
\caption{The real part of the situation in the Conrey-Ghosh case. The circle at
  $(1,1/2)$ indicates the point for which the value of
  $\mathcal{J}(2k,2h)$ is coveted. The dots indicate the
locations where Expression~(\ref{computational}) applies, and the
crosses indicate the points to which that expression is
meromorphically continued (for a fixed $h$, i.e. horizontally) thanks
to Expression~(\ref{AtOne}). Note that for fixed integer $h$, 
this continuation hits a pole when crossing the dashed line (and many
more before reaching $k=1$, as $h$ increases: see Figure~\ref{fig_one}).}
\label{ConreyGhoshFig}
\end{figure}

We offer in Figure~\ref{ConreyGhoshFig} one way to circumvent this
problem. The idea is to compute for each fixed integer $h$ the values
of the meromorphic continuation in $k$ of $(\mathcal{M})(2k,2h)$ at
$k=1$ (i.e. at the crosses). This should be enough to know through
Equation~(\ref{relationTwo}) any value of the form
$\left|\mathcal{M}\right|(2,2h)$, which could then finally be used to 
meromorphically continue  $\left|\mathcal{V}\right|(2,2h)$ to $h=1/2$.

Getting the meromorphic continuation of Equation~(\ref{computational})
to $k=1$ is quite subtle. 
\begin{prop} Define $ (\mathcal{M})(2,r)$ as the meromorphic
  continuation in $k$ of $ (\mathcal{M})(2k,r)$, evaluated at
  $k=1$. Then, the exponential generating series of $
  (\mathcal{M})(2,r)$ is given by
\label{propAtOne}
\begin{multline}
\label{AtOne}
\sum_{r>0} (\mathcal{M})(2,r) \frac{(\mathfrak{i}z)^r}{r!}
=
\\
\sum_{d=1}^\infty
\sum_{\substack{\sort{s},\sort{t} \in \mathbb{N}^d}}
\left|\frac{ 1 }{ s_i!  t_j!(1+s_i+t_j)}\right|_{d\times d}^2 
\left(\prod_{i,j=1}^d
\frac{1}{2^{\nu(t_j)}}\frac{1-t_j}{2+s_j}
\right)
z^{d+\total s_i+t_i},
\end{multline}
where $\nu(0)=0$ when $t=0$, $\nu(t)=1$ when $t\ge 2$. The value
$\nu(1)$ is free to choose.
\begin{proof}
When looking for the analytic continuation in $k$, most of the formulas we
have found so far are misleading. For instance, in light of the remark
in footnote~\ref{length}, one could think that the sums over
partitions $\lambda$ in Expression~(\ref{SchurAppear}) or
Theorem~\ref{ninthTheorem} immediately reduce when $k=1$ to sums over
partitions $\lambda$ of length 1, i.e. partitions indexed by a single
variable. However, in those cases, the other factor in the summands (i.e. for
instance $\tilde{s}_\lambda^{(2k)}$ in Theorem~\ref{ninthTheorem})
might actually be undefined if we take $k=1$ (in that particular case, when
$l(\lambda)\ge 3$\footnote{or equivalently when $t_1\ge2$ if
  $\lambda=\Frob{\sort{s}}{\sort{t}}$}). 

 We can get a better intuition through
 Expression~(\ref{computational}), which we use as a basis of our
 proof. We are clearly required to find the  meromorphic continuation
to $k=1$ and for fixed $s,t\ge 0\in \mathbb{N}$
of $$\frac{\Gamma(2k-t)}{\Gamma(k-t)}\cdot\frac{\Gamma(k+s+1)}{\Gamma(2k+s+1)}.$$
The second factor is certainly not a problem and immediately gives
$\frac{1}{s+2}$. For the first factor, we have to look at 
$\lim_{k \rightarrow 1}\frac{\Gamma(2k-t)}{\Gamma(k-t)}$
for $t\ge 0$. Pick any integer $a$ such that $1+a-t \ge 0$. Then,
using the functional equation for $\Gamma$, we have
$$
\lim_{k \rightarrow 1}\frac{\Gamma(2k-t)}{\Gamma(k-t)} = 
\lim_{k \rightarrow 1}\frac{\Gamma(2k+a-t)}{\Gamma(k+a-t+1)} \cdot \frac{(k-t)(k-t+1)\cdots(k-1)\cdots(k+a-t)}{(2k-t)\cdots(2k-2)\cdots(2k+a-t-1)}.
$$
Note that the terms $\frac{k-1}{2k-2}$ only appear if $t\ge 2$. In
that case we get
$$
\lim_{k \rightarrow 1}\frac{\Gamma(2k-t)}{\Gamma(k-t)} = 
(1-t) \lim_{k \rightarrow 1}\frac{\Gamma(2k+a-t)}{\Gamma(k+a-t+1)} \cdot
\frac{k-1}{2k-2}=\frac{1}{2} (1-t), 
$$
and in the case $t \le 2$ the factor of 2 is missing.
\end{proof}
\end{prop}

One can also check that the values recovered using
Proposition~\ref{propAtOne} agree with the values obtained using
$\frac{X_r(2)}{Y_r(2)}$ and thus Theorem~\ref{thmComplexM}.

For completeness, we give the beginning of the sequence of $X_r(2)$s,
for $r=1$ to~$15$:
\begin{multline*}
    1,
    2,
    0,
    18,
    50,
    -6540,
    -11760,
    852180,
    1228500,
    590126040,
    558613440,\\
    -39273224760,
    455842787400,
    5775116644337040,
    14904865051876800
\end{multline*}

Unfortunately, we fall short of actually finding the full meromorphic
continuation of $(\mathcal{M})(2,r)$ and have to leave this for a
further paper. 
\section{Conclusion}
The initial goal was to compute the $\left(\mathcal{M}\right)(2k,r)$,
$\left|\mathcal{M}\right|(2k,2h)$ and $|\mathcal{V}|(2k,2h)$ more
effectively than previously done. 

We feel that we have achieved this goal, since we have been able to
shed some light (for instance in Theorem~\ref{thmComplexM}) on the
structure of the results. This structure (rational functions with
known denominators)  underlines tables already  available in
\cite{HughesJoint} or \cite{CRS}. We have also been able to use these
results to obtain better algorithms to compute those rational
functions, thereby extending the data that was available. \ForArxiv{Much
  of that data is now available in the source of the
  \arXiv\,submission, or at \cite{DehayeData}.} As a
corollary we have shown that for large(r) $h$ the roots (in $k$) of
$|\mathcal{V}|(2k,2h)$ cease to all be real, a fluke only for the small-$h$
cases available previously.

However, we have not obtained a formula for \emph{all}
$|\mathcal{V}|(2k,r)$. In particular, we cannot recover the value of
$|\mathcal{V}|(2,1)$, which can be conjectured from Conrey and Ghosh's
result for $\mathcal{J}(2,1)$. 

Those methods should also give more general moments, for instance for
expressions of the form 
\begin{eqnarray*}
\average{\left|Z_U(\theta_1)\right|^{2k} \left|
    \frac{Z_U'(\theta_2)}{Z_U(\theta_2)} \right|^{r}}
\end{eqnarray*}
or
\begin{eqnarray}
\average{\left|Z_U(\theta_1)\right|^{2k} \left|
    \frac{Z_U''(\theta_2)}{Z_U(\theta_2)} \right|^{r}}.
\label{secondderiv}
\end{eqnarray}
An expression for those two extensions in the shape of Equation~(\ref{RealMomentSeries}) would
definitely be available (for instance, in the case of
Expression~(\ref{secondderiv}), we would most likely have to compute the
equivalent of Equation~(\ref{RealMomentSeries}) by summing over
$\sort{\boldmu} \in (2\Nplus)^r$). However, the second part of the
computation, the part covered here by
Proposition~\ref{propStructuresCount}, would probably be significantly
worsened. 
\bibliographystyle{alpha}
\bibliography{references}
\end{document}